# A FREQUENCY DOMAIN EMPIRICAL LIKELIHOOD FOR SHORT- AND LONG-RANGE DEPENDENCE

By Daniel J. Nordman and Soumendra N. Lahiri

*Iowa State University*

This paper introduces a version of empirical likelihood based on the periodogram and spectral estimating equations. This formulation handles dependent data through a data transformation (i.e., a Fourier transform) and is developed in terms of the spectral distribution rather than a time domain probability distribution. The asymptotic properties of frequency domain empirical likelihood are studied for linear time processes exhibiting both short- and long-range dependence. The method results in likelihood ratios which can be used to build nonparametric, asymptotically correct confidence regions for a class of normalized (or ratio) spectral parameters, including autocorrelations. Maximum empirical likelihood estimators are possible, as well as tests of spectral moment conditions. The methodology can be applied to several inference problems such as Whittle estimation and goodness-of-fit testing.

**1. Introduction.** The main contribution of this paper is a new formulation of empirical likelihood (EL) for inference with two fundamentally different types of dependent data: time series exhibiting either short-range dependence (SRD) or long-range dependence (LRD). Let $\{X_t\}, t \in \mathbb{Z}$, be a stationary sequence of random variables with mean $\mu$ and spectral density $f$ on $\Pi = [-\pi, \pi]$, where

$$(1) \qquad f(\lambda) \sim C(\alpha)|\lambda|^{-\alpha}, \qquad \lambda \to 0,$$

for $\alpha \in [0, 1)$ and a constant $C(\alpha) > 0$ involving $\alpha$ (with $\sim$ indicating that the terms have a ratio of one in the limit). When $\alpha = 0$, we classify the process $\{X_t\}$ as *short-range dependent* (SRD). For $\alpha > 0$, the process will be called *long-range dependent* (LRD). This classification resembles that of [25] and encompasses the formulation of LRD described in [4, 40].

Originally proposed by [33, 34] for independent samples, EL allows for nonparametric likelihood-based inference with a broad range of applications [36]. An important benefit of EL inference is that confidence regions for parameters may be calibrated through log-likelihood ratios, without requiring any direct estimates of variance or skewness [19]. However, a difficulty with extending EL methods to dependent data is then to ensure that "correct" variance estimation occurs automatically within EL ratios under the data dependence structure. This is an important reason why the EL version for i.i.d. data from [34] does not apply to dependent data (see [23], page 2085).

Recent extensions of EL to time series in [23, 31] have relied exclusively on a SRD structure with rapidly decreasing process correlations. In particular, [23] provided a breakthrough formulation of EL for weakly dependent data based on data blocks rather than individual observations. Under SRD, the resulting blockwise EL ratios correctly perform variance estimation of sample means within their mechanics. Data blocking has also proven to be crucial in extending other nonparametric likelihoods to weakly dependent processes, such as the block bootstrap and resampling methods described in [26], Chapter 2.

In comparison to weak dependence, the rate of decay of the covariance function $r(k) = \text{Cov}(X_j, X_{j+k})$ is characteristically much slower under strong dependence $\alpha > 0$, namely,

$$(2) \qquad r(k) \sim \tilde{C}(\alpha) k^{-(1-\alpha)}, \qquad k \to \infty,$$

with a constant $\tilde{C}(\alpha) > 0$, which is an alternative representation of LRD (with equivalence to (1) if the covariances converge quasimonotonically to zero; see page 1632 of [40]). This autocovariance behavior implies that statistical procedures developed for SRD may not be applicable under LRD, often due to complications with variance estimation. For example, the moving block bootstrap is known to be invalid under strong dependence for inference on the process mean $EX_t = \mu$ [24], partly because (2) implies that the variance $\text{Var}(\bar{X}_n) = O(n^{-1+\alpha})$ of a size $n$ sample mean exhibits a slower, *unknown* rate of decay compared to the $O(n^{-1})$ rate associated with SRD data. For this reason, the blockwise EL formulation of [23] will also break down under strong dependence for inference on the mean.

In this paper, we formulate an EL based on the periodogram combined with certain estimating equations. Using this data transformation to weaken the underlying dependence structure, the resulting frequency domain empirical likelihood (FDEL) provides a common tool for nonparametric inference on both SRD *and* LRD time series. Because this EL version involves the spectral distribution of a time process rather than a time domain probability distribution, inference is restricted to a class of normalized spectral parameters described in Section 2. The frequency domain bootstrap (FDB) of [11], developed for SRD, targets the same class of parameters. Hence, for parameters *not* defined in terms of the spectral density (e.g., the process mean $\mu$), the FDEL is inapplicable, while the time domain blockwise EL in [23] may still be valid if the process exhibits SRD (i.e., is valid for a larger class of parameters under *weak* dependence).

Our main result is the asymptotic distribution of FDEL ratio statistics, which are shown to have limiting chi-square distributions under both SRD and LRD for setting confidence regions. That is, FDEL shares the strength of EL methods to incorporate "correct" variance estimation for spectral parameter inference automatically in its mechanics. For normalized spectral parameters where both the FDB and blockwise EL may be applicable under SRD (e.g., autocorrelations), the FDEL requires no kernel density estimates of $f$ (as with the FDB) and no block selection (as with the blockwise EL). Our FDEL results also refine some EL theory given in [31], wherein periodogram-based EL confidence regions for Whittle-type estimation with SRD linear processes were proposed. We additionally consider FDEL tests based on maximum EL estimation, which are helpful for assessing both parameter conjectures and the validity of (spectral) moment conditions, as in the independent data EL formulation [34, 38].

The methodology presented here is applicable to linear processes with spectral densities satisfying (1), which includes two common models for LRD processes: the fractional Gaussian processes of [29] with spectral density

$$
\begin{aligned}
(3) \qquad f_{H,\sigma^2}(\lambda) &= \frac{4\sigma^2\Gamma(2H-1)}{(2\pi)^{2H+2}}\cos(\pi H - \pi/2)\sin^2(\lambda/2)\\
&\quad \times \sum_{k=-\infty}^{\infty} |\lambda/(2\pi)+k|^{-1-2H}, \qquad \lambda \in \Pi,
\end{aligned}
$$

$1/2 < H < 1$, and the fractional autoregressive integrated moving average (FARIMA) processes of [1, 18, 21] with spectral density

$$
(4) \qquad f_{d,\rho,\varrho,\sigma^2}(\lambda) = \frac{\sigma^2}{2\pi}|1-e^{\imath\lambda}|^{-2d}\left|\frac{\sum_{j=0}^{p}\rho_j(e^{\imath\lambda})^j}{\sum_{j=0}^{q}\varrho_j(e^{\imath\lambda})^j}\right|^2, \qquad \lambda \in \Pi,
$$

based on parameters $0 < d < 1/2$, $\rho = (\rho_1,\ldots,\rho_p)$, $\varrho = (\varrho_1,\ldots,\varrho_q)$ with $\rho_0 = \varrho_0 = 1$ and $\imath = \sqrt{-1}$. These models fulfill (1) with $\alpha = 2H - 1$ and $\alpha = 2d$, respectively.

The rest of the paper is organized as follows. Section 2 describes the role of spectral estimating equations for FDEL inference and provides several examples. In Section 3, we explain the construction of EL in the frequency domain. Section 4 contains the assumptions and the main results on the distribution of FDEL log-ratios for confidence region estimation and simple hypothesis testing. In Section 5, we consider maximum EL estimation in the frequency domain. We describe the application of FDEL to Whittle estimation in Section 6, while Section 7 considers goodness-of-fit testing with FDEL. Section 8 offers some conclusions. Proofs of the results are given in Section 9 and the Appendix.

**2. Spectral estimating equations.** Consider inference on a parameter $\theta \in \Theta \subset \mathbb{R}^p$ based on a time stretch $X_1, \ldots, X_n$. Following the EL framework of [38, 39] with i.i.d. data, we suppose that information about $\theta$ exists through a system of general estimating equations. However, we will use the process *spectral distribution* to define moment conditions as follows. Let

$$(5) \qquad G_\theta(\lambda) = (g_{1,\theta}(\lambda), \ldots, g_{r,\theta}(\lambda))' : \Pi \times \Theta \to \mathbb{R}^r$$

denote a vector of even, estimating functions with $r \geq p$. For the case $r > p$, the above functions are said to be "overidentifying" for $\theta$. We assume that $G_\theta$ satisfies the spectral moment condition

$$(6) \qquad \int_0^\pi G_{\theta_0}(\lambda) f(\lambda) \, d\lambda = \mathcal{M}$$

for some known $\mathcal{M} \in \mathbb{R}^r$ at the true value $\theta_0$ of the parameter. As distributional results in Section 4 indicate, we will typically require $\mathcal{M} = 0$, which places some restrictions on the types of spectral parameters considered. However, the FDEL framework is valid for estimating normalized spectral means: $\theta = \int_0^\pi G f \, d\lambda / \int_0^\pi f \, d\lambda$, based on a vector function $G$. The FDB targets the same parameters under SRD; [11] comments on the importance, and often complete adequacy, of population information expressed in this ratio form. The FDEL construction in Section 3 combines the periodogram with the estimating equations in (6).

2.1. *Examples.* We provide a few examples of useful estimating functions for inference, some of which satisfy (6) with $\mathcal{M} = 0$.

EXAMPLE 1 (Autocorrelations). Consider interest in the autocorrelation function $\rho(\cdot)$ at arbitrary lags $m_1, \ldots, m_p$, that is, $\theta = (\rho(m_1), \ldots, \rho(m_p))'$, where

$$\rho(m) = r(m)/r(0) = \int_0^\pi \cos(m\lambda) f(\lambda) \, d\lambda \bigg/ \int_0^\pi f(\lambda) \, d\lambda, \qquad m \in \mathbb{Z}.$$

One can select $G_\theta(\lambda) = (\cos(m_1\lambda), \ldots, \cos(m_p\lambda))' - \theta$ for autocorrelation inference, fulfilling (6) with $\mathcal{M} = 0 \in \mathbb{R}^p$ and $r = p$.

EXAMPLE 2 (Spectral distribution function). For $\omega \in [0, \pi]$, denote the spectral distribution function by $F(\omega) = \int_0^\omega f(\lambda) \, d\lambda$. Suppose $\theta = (F(\tau_1)/F(\pi), \ldots, F(\tau_p)/F(\pi))'$ for some $\tau_1, \ldots \tau_p \in (0, \pi)$. This normalized parameter $\theta$ often sufficiently characterizes the spectral distribution $F$ for testing purposes [9]. For inference on $\theta$, we can pick $G_\theta(\lambda) = (\mathbb{1}\{\lambda \leq \tau_1\}, \ldots, \mathbb{1}\{\lambda \leq \tau_p\})' - \theta$ where $\mathbb{1}\{\cdot\}$ denotes the indicator function. Then (6) holds with spectral mean $\mathcal{M} = 0 \in \mathbb{R}^p$.

EXAMPLE 3 (Goodness-of-fit tests). There has been increasing interest in frequency-domain-based tests to assess model adequacy [2, 37]. Consider a test involving a simple null hypothesis $H_0 : f = f_0$ against an alternative $H_1 : f \neq f_0$ for some candidate density $f_0$. With EL techniques, one immediate test for $H_0$ is based on the function $G_0(\lambda) = 1/f_0(\lambda)$ with spectral mean $\pi$ under $H_0$; here, we treat $r = 1$ and the dimension $p$ of $\theta$ as 0. We show in Section 7 that this results in an EL ratio test which resembles a spectral goodness-of-fit test statistic proposed by [30] and shown by [3] to be useful for LRD Gaussian series. The more interesting and complicated problem of testing the hypothesis that $f$ belongs to a given model family can also be addressed with FDEL, as discussed in Section 7.

EXAMPLE 4 (Whittle estimation). We denote a parametric collection of spectral densities by

$$\mathcal{F} = \{f_\theta(\lambda) : \theta \in \Theta\} \tag{7}$$

and assume the densities are positive on $\Pi$ and identifiable [e.g., for $\theta \neq \tilde{\theta} \in \Theta$, the Lebesgue measure of $\{\lambda : f_\theta(\lambda) \neq f_{\tilde{\theta}}(\lambda)\}$ is positive]. For fitting the model $f_\theta$ to the data, Whittle estimation [42] seeks the $\theta$-value at which the theoretical "distance" measure

$$W(\theta) = (4\pi)^{-1} \int_0^\pi \left\{ \log f_\theta(\lambda) + \frac{f(\lambda)}{f_\theta(\lambda)} \right\} d\lambda \tag{8}$$

achieves its minimum [13]. The model class may be misspecified (possibly $f \notin \mathcal{F}$), but Whittle estimation aims for the density in $\mathcal{F}$ "closest" to $f$, as measured by $W(\theta)$.

To consider a particular parameterization of (7), suppose

$$\begin{aligned} f_\theta(\lambda) &= \sigma^2 k_\vartheta(\lambda), \\ \theta &= (\sigma^2, \vartheta')', \Theta \subset (0, \infty) \times \mathbb{R}^{p-1}, \ \vartheta = (\vartheta_1, \ldots, \vartheta_{p-1})', \end{aligned} \tag{9}$$

with kernel density $k_\vartheta$, and that Kolmogorov's formula holds with $(2\pi)^{-1} \times \int_{-\pi}^\pi \log f_\theta(\lambda) \, d\lambda = \log[\sigma^2/(2\pi)]$ (e.g., taking $\sigma^2$ as the innovation variance in a linear model). The model class in (9) is commonly considered in the context of Whittle estimation for both SRD and LRD time series, including those LRD processes formulated in (3) and (4) (see [10, 15, 17, 20]). Under appropriate conditions, the true minimum argument $\theta_0 = (\sigma_0^2, \vartheta_0')'$ of $W(\theta)$ is determined by the stationary solution of $\partial W(\theta)/\partial \theta = 0$ or

$$\int_0^\pi f(\lambda)\{\partial f_\theta^{-1}(\lambda)/\partial \vartheta\} \, d\lambda = 0, \qquad \pi^{-1} \int_0^\pi f(\lambda) f_\theta^{-1}(\lambda) \, d\lambda = 1, \tag{10}$$

where $f_\theta^{-1}(\lambda) \equiv 1/f_\theta(\lambda)$. The moment conditions in (10) give a set of estimating functions for FDEL inference on $\theta$ defining the densities in (9). Namely, the choice

$$
\begin{aligned}
(11) \quad & G_\theta^w(\lambda) = (f_\theta^{-1}(\lambda), \partial f_\theta^{-1}(\lambda)/\partial \vartheta_1, \ldots, \partial f_\theta^{-1}(\lambda)/\partial \vartheta_{p-1})', \\
& \mathcal{M}_w = (\pi, 0, \ldots, 0)' \in \mathbb{R}^p,
\end{aligned}
$$

fulfills (6). The FDB uses similar estimating equations for Whittle parameter inference [11]. To treat $\sigma^2$ as a nuisance parameter, which is common for densities as in (9), estimating functions

$$
(12) \quad G_\vartheta^{w*}(\lambda) = \partial k_\vartheta^{-1}(\lambda)/\partial \vartheta, \qquad \mathcal{M}_{w*} = 0 \in \mathbb{R}^{p-1},
$$

provide structure for inference on the remaining parameters $\vartheta$ determining $k_\vartheta$ in $f_\theta$.

**3. Definition of frequency domain empirical likelihood.** Denote the periodogram of the sequence $X_1, \ldots, X_n$ by $I_n(\lambda) = (2\pi n)^{-1} |\sum_{t=1}^n X_t \exp(-\imath t \lambda)|^2$, $\lambda \in \Pi = [-\pi, \pi]$, where $\imath = \sqrt{-1}$. Using estimating functions $G_\theta$ as in (5), the profile FDEL function for $\theta \in \Theta$ is given by

$$
(13) \quad L_n(\theta) = \sup\left\{ \prod_{j=1}^N w_j : w_j \geq 0, \; \sum_{j=1}^N w_j = \pi, \; \sum_{j=1}^N w_j G_\theta(\lambda_j) I_n(\lambda_j) = \mathcal{M} \right\},
$$

where $\lambda_j = 2\pi j/n$, $j \in \mathbb{Z}$, are Fourier frequencies and $N = \lfloor (n-1)/2 \rfloor$. Point masses $w_j$ assigned to each ordinate $\lambda_j$ create a discrete measure on $[0, \pi]$ with the restriction that the integral of $G_\theta I_n$, based on this measure, equals $\mathcal{M}$. The largest possible product of these point masses determines the FDEL function for $\theta \in \Theta$. When the conditioning set in (13) is empty, we define $L_n(\theta) = -\infty$. If $\mathcal{M}$ is interior to the convex hull of $\{\pi G_\theta(\lambda_j) I_n(\lambda_j)\}_{j=1}^N$, then $L_n(\theta)$ is a positive constrained maximum solved by optimizing

$$
\mathcal{L}(w_1, \ldots, w_N, \gamma, t) = \sum_{j=1}^N \log(w_j) + \gamma\left(\pi - \sum_{j=1}^N w_j\right) \\
- Nt'\left(\sum_{j=1}^N w_j G_\theta(\lambda_j) I_n(\lambda_j) - \mathcal{M}\right),
$$

with Lagrange multipliers $\gamma$ and $t = (t_1, \ldots, t_r)'$ as in [33, 34]. Then (13) may be written as

$$
\begin{aligned}
(14) \quad & L_n(\theta) = \pi^N \prod_{j=1}^N p_j(\theta), \\
& p_j(\theta) = N^{-1}[1 + t_\theta'\{\pi G_\theta(\lambda_j) I_n(\lambda_j) - \mathcal{M}\}]^{-1} \in (0, 1),
\end{aligned}
$$

where $t_\theta$ is the stationary point of the function $q(t) = \sum_{j=1}^{N} \log(1 + t'\{\pi G_\theta(\lambda_j) \times I_n(\lambda_j) - \mathcal{M}\})$. (See [34, 38] for further computational details on EL.) Without the integral-type linear constraint in (13), $\prod_{j=1}^{N} w_j$ has a maximum when each $w_j = \pi/N$, so we can form a profile EL ratio

$$(15) \quad R_n(\theta) = L_n(\theta)/(\pi N^{-1})^N = \prod_{j=1}^{N}[1 + t'_\theta\{\pi G_\theta(\lambda_j)I_n(\lambda_j) - \mathcal{M}\}]^{-1}.$$

3.1. *A density-based formulation of empirical likelihood.* To help relate the EL results here to those in [31], we give an alternative, model-based formulation of FDEL. This version requires a density class $\mathcal{F}$ as in (7) and involves approximating the expected value $\mathrm{E}(I_n(\lambda_j))$ with $f_\theta(\lambda_j)$, using a density $f_\theta \in \mathcal{F}$. Namely, let

$$(16) \quad \begin{aligned} L_{n,\mathcal{F}}(\theta) &= \sup\Bigg\{\prod_{j=1}^{N} w_j : w_j \geq 0, \ \sum_{j=1}^{N} w_j = \pi, \\ &\qquad \sum_{j=1}^{N} w_j G_\theta(\lambda_j)[I_n(\lambda_j) - f_\theta(\lambda_j)] = 0\Bigg\} \\ R_{n,\mathcal{F}}(\theta) &= (N/\pi)^N L_{n,\mathcal{F}}(\theta). \end{aligned}$$

We consider the densities $f_\theta$ and prospective functions $G_\theta$ as dependent on the same parameters, which causes no loss of generality. An exact form for $L_{n,\mathcal{F}}(\theta)$ can be deduced as with $L_n(\theta)$, obtained by replacing $I_n(\lambda_j), \mathcal{M}$ with $I_n(\lambda_j) - f_\theta, 0$ in (14).

Section 6 discusses the model-based EL ratio in (16) for refining results in [31] on confidence interval estimation of Whittle parameters. Additionally, this version of FDEL may be suitable for conducting goodness-of-fit tests with respect to a family of spectral densities.

**4. Main result: distribution of empirical likelihood ratio.** Before describing the distributional properties of FDEL, we provide some assumptions on the time process under consideration and the potential vector of estimating functions $G_\theta$ in (5).

4.1. *Assumptions.* In the following, let $\theta_0$ denote the unique (true) parameter value which satisfies (6).

ASSUMPTION A.1. $\{X_t\}$ is a real-valued, linear process with a moving average representation of the form

$$X_t = \mu + \sum_{j=-\infty}^{\infty} b_j \varepsilon_{t-j}, \qquad t \in \mathbb{Z},$$

where $\{\varepsilon_t\}$ are i.i.d. random variables with $\mathrm{E}(\varepsilon_t) = 0$, $\mathrm{E}(\varepsilon_t^2) = \sigma_\varepsilon^2 > 0$, $\mathrm{E}(\varepsilon_t^8) < \infty$ and fourth order cumulant denoted by $\kappa_{4,\varepsilon} \equiv \mathrm{E}(\varepsilon_t^4) - 3\sigma_\varepsilon^4$, $\{b_t\}$ is a sequence of constants satisfying $\sum_{t \in \mathbb{Z}} b_t^2 < \infty$ and $b_0 = 1$, and $f(\lambda) = \sigma_\varepsilon^2 |b(\lambda)|^2/(2\pi)$, $\lambda \in \Pi$, with $b(\lambda) = \sum_{j \in \mathbb{Z}} b_j e^{\iota j \lambda}$. It is assumed that $f(\lambda)$ is continuous on $(0, \pi]$ and that $f(\lambda) \leq C|\lambda|^{-\alpha}$, $\lambda \in \Pi$, for some $\alpha \in [0, 1)$, $C > 0$.

ASSUMPTION A.2. Each component $g_{j,\theta_0}$ of $G_{\theta_0}$ is an even, integrable function such that $|g_{j,\theta_0}(\lambda)| \leq C|\lambda|^\beta$, $\lambda \in \Pi$, where $0 \leq \beta < 1$, $\alpha - \beta < 1/2$, $j = 1, \ldots, r$.

ASSUMPTION A.3. For each $g_{j,\theta_0}(\lambda)$, $j = 1, \ldots, r$, one of the following is fulfilled:

CONDITION 1. $g_{j,\theta_0}$ is Lipschitz of order greater than $1/2$ on $[0, \pi]$.
CONDITION 2. $g_{j,\theta_0}$ is continuous on $\Pi$ and $|\partial g_{j,\theta_0}(\lambda)/\partial \lambda| \leq C|\lambda|^{\beta_j - 1}$ for some $0 \leq \beta_j < 1$, $2\alpha - \beta_j < 1$.
CONDITION 3. $g_{j,\theta_0}$ is of bounded variation on $[0, \pi]$ with finite discontinuities and $\alpha < 1/2$, with $|r(k)| \leq Ck^{-\upsilon}$ for some $\upsilon > 1/2$ (e.g., $\upsilon = 1 - \alpha$).

ASSUMPTION A.4. The $r \times r$ matrix $W_{\theta_0} = \int_\Pi f^2(\lambda) G_{\theta_0}(\lambda) G'_{\theta_0}(\lambda) \, d\lambda$ is positive definite.

ASSUMPTION A.5. On $(0, \pi]$, either (i) $f$ is differentiable and $|\partial f(\lambda)/\partial \lambda| \leq C|\lambda|^{-\alpha-1}$ or (ii) each $f(\lambda)g_{j,\theta_0}(\lambda)$ is of bounded variation or is piecewise Lipschitz of order greater than $1/2$ on $[0, \pi]$, $j = 1, \ldots, r$. As $n \to \infty$, $P(0 \in \mathrm{ch}^\circ \{\pi G_{\theta_0}(\lambda_j)[I_n(\lambda_j) - f(\lambda_j)]\}_{j=1}^N) \to 1$, where $\mathrm{ch}^\circ A$ denotes the interior convex hull of a finite set $A \subset \mathbb{R}^r$.

We briefly discuss the assumptions. The bound on $f$ in Assumption A.1 allows for the process $\{X_t\}$ to exhibit both SRD and LRD and is a slight generalization of (1). The behavior of $G_{\theta_0}$ in Assumption A.2 controls the growth rate of the scaled periodogram ordinates, $G_{\theta_0}(\lambda_j) I_n(\lambda_j)$, at low frequencies under LRD and ensures that $W_{\theta_0}$ is finite. Important processes are permissible under A.1 and for these, useful estimating functions often satisfy A.2. Assumption A.3 outlines smoothness criteria for the estimating functions. The estimating functions treated by the FDB in [11] satisfy A.3, including those for autocorrelations and normalized spectral distribution in Section 2.1. The functions $f_\theta^{-1}$ and $\partial f_\theta^{-1}/\partial \theta$ from Examples 3 and 4 satisfy A.3 for use in Whittle-like estimation and goodness-of-fit testing with many SRD and LRD models in the FDEL framework. For example, Hannan [20] considers Whittle estimation for ARMA densities for which functions $G_\theta^w$ in (11) satisfy Condition 1. The functions $f_\theta^{-1}$ and $\partial f_\theta^{-1}/\partial \theta$ associated with the fractional Gaussian and FARIMA LRD densities in (3) and (4) fulfill

Condition 2 [10, 15, 17]. Process dependence that is not extremely strong, so that $f^2$ is integrable, allows greater flexibility in choosing more general estimating functions in Condition 3.

For EL inference exclusively with the model-based functions $L_{n,\mathcal{F}}$ or $R_{n,\mathcal{F}}$ from Section 3.1, we introduce the additional assumption A.5 which is generally not restrictive. The probabilistic condition in A.5 implies only that the EL ratio $R_{n,\mathcal{F}}$ can be finitely computed at $\theta_0$, resembling EL assumptions from [31] and [35].

4.2. *Asymptotic distribution of empirical likelihood ratio and confidence regions.* We now establish a nonparametric recasting of Wilks' theorem [43] for FDEL ratios under SRD and LRD, useful for setting confidence regions and making simple hypothesis tests, as in [33, 34, 35]. Define two scaled log-profile FDEL ratio statistics

$$(17) \qquad \ell_n(\theta) = -4 \log R_n(\theta) \quad \text{and} \quad \ell_{n,\mathcal{F}}(\theta) = -2 \log R_{n,\mathcal{F}}(\theta),$$

using (15) and (16). The difference in the scalar adjustments to log-likelihoods in (17) is due to the assumption that the periodogram ordinates are "mean-corrected" in the construction of $\ell_{n,\mathcal{F}}(\theta)$. In the following theorem, $\chi^2_\nu$ denotes a chi-square distribution with $\nu$ degrees of freedom:

THEOREM 1. *Suppose Assumptions A.1–A.4 hold. If $\mathcal{M} = 0 \in \mathbb{R}^r$, then as $n \to \infty$,*

(i) $\ell_n(\theta_0) \xrightarrow{d} \chi^2_r$.
(ii) *Additionally, if A.5 holds and $f = f_{\theta_0}$, then $\ell_{n,\mathcal{F}}(\theta_0) \xrightarrow{d} \chi^2_r$.*
(iii) *If $\kappa_{4,\varepsilon} = 0$, statement* (ii) *remains valid, even if $\mathcal{M} \neq 0 \in \mathbb{R}^r$.*

REMARK 1. For a Gaussian $\{X_t\}$ process, the fourth order innovation cumulant $\kappa_{4,\varepsilon} = 0$.

Due to the data transformation aimed at weakening the time dependence structure, FDEL ratios closely resemble EL ratios with i.i.d. data [34, 38]. The formulation of estimating equations satisfying $\mathcal{M} = 0$ in (6) is generally necessary for $\ell_n(\theta_0)$ to have a chi-square limit and is a consequence of this EL based on the periodogram. A similar moment restriction is shared by the FDB, as detailed in [11] (page 1938), due to difficulties in estimating the variance of empirical spectral means. Similar complications arise in the inner mechanics of FDEL, requiring $\mathcal{M} = 0$. As the proof of Theorem 1 shows [see (24) and (27)], variance estimators intrinsic to FDEL ratios are asymptotically of the form given in Lemma 7 of Section 9 (i.e., setting $gh = GG'$ there) and these target the asymptotic variance $V$ of an empirical spectral mean appearing in Lemma 6 so that variance estimation within FDEL is consistent if $\mathcal{M} = 0$ (or if the innovation cumulant $\kappa_{4,\varepsilon} = 0$); see [32] for details. However, Section 2.1 gives some important estimating equations for which $\mathcal{M} = 0$ and, importantly, estimating functions may be chosen with more flexibility for inference on Gaussian processes.

If the true density $f$ belongs to $\mathcal{F}$ in (7), then a confidence region can be calibrated with $\ell_{n,\mathcal{F}}(\theta)$ as well. Monti [31] suggests similar confidence regions for Whittle estimation with SRD linear processes and spectral densities parameterized as in (9) (e.g., ARMA models). However, if the candidate density class $\mathcal{F}$ is incorrect, confidence regions based on $\ell_{n,\mathcal{F}}(\theta)$ become conservative to a degree dependent on the misspecification. This closely parallels the behavior of the EL ratio with misspecified regression models, as described in Section 5.4 of [35]. Confidence regions set with $\ell_n$ do not generally require specification of a model density class $\mathcal{F}$, but for the case of inference on Whittle parameters where a class of densities $\mathcal{F}$ may be involved, Section 6 describes how $\ell_n$ (unlike $\ell_{n,\mathcal{F}}$) may be used, even if $\mathcal{F}$ is misspecified.

**5. Extensions to maximum empirical likelihood estimation.** We shall refer to the maximum of $R_n(\theta)$ from (15) as the *maximum empirical likelihood estimator* (MELE) and denote it by $\hat{\theta}_n$; we denote the maximum of $R_{n,\mathcal{F}}(\theta)$ from (16) by $\hat{\theta}_{n,\mathcal{F}}$. We next show that with both SRD and LRD linear time series, maximum empirical likelihood estimates (MELE's) $\hat{\theta}_n$ and $\hat{\theta}_{n,\mathcal{F}}$ have properties resembling those available in EL frameworks involving independent data.

5.1. *Consistency and asymptotic normality.* We first consider establishing the existence, consistency and asymptotic normality of a sequence of local maximums of FDEL functions $R_n(\theta)$ and $R_{n,\mathcal{F}}(\theta)$, along the lines of the classical arguments of [7]. The assumptions involved are very mild and have the advantage of being typically simple to verify and apply to both versions of FDEL $R_n(\theta)$ and $R_{n,\mathcal{F}}(\theta)$. Qin and Lawless [38] adopt a similar approach to study MELE's in the i.i.d. data scenario.

Let $\|\cdot\|$ denote the Euclidean norm. For $n \in \mathbb{N}$, define the open neighborhood $\mathcal{B}_n = \{\theta \in \Theta : \|\theta - \theta_0\| < n^{-\eta}\}$, where $\eta = \max\{1/3, 1/4 + (\alpha - \beta)/2, (1 + \alpha + \delta)/4\} < 1/2$ for $\delta < 1$ defined below in Theorem 2 and $\alpha$, $\beta$ from Assumptions A.1–A.2.

THEOREM 2. *Assume* A.1–A.4 *hold and* $\mathcal{M} = 0$. *Suppose, in a neighborhood of* $\theta_0$, $\partial G_\theta(\lambda)/\partial \theta$, $\partial^2 G_\theta(\lambda)/\partial \theta \partial \theta'$ *are continuous in* $\theta$ *and* $\|\partial G_\theta(\lambda)/\partial \theta\|$, $\|\partial^2 G_\theta(\lambda)/\partial \theta \partial \theta'\|$ *are bounded by* $C|\lambda|^{-\delta}$ *for some* $\delta < 1$, $\delta + \alpha < 1$. *Suppose further that* $\partial G_{\theta_0}/\partial \theta$ *is Riemann integrable and that* $D_{\theta_0} \equiv \int_\Pi f(\lambda) \partial G_{\theta_0}(\lambda)/\partial \theta d\lambda$ *has full column rank* $p$.

(i) *As $n \to \infty$, there exists a sequence of statistics $\{\hat{\theta}_n\}$ such that $P$ ($\hat{\theta}_n$ is a maximum of $R_n(\theta)$ and $\hat{\theta}_n \in \mathcal{B}_n$) $\to 1$ and*

$$\sqrt{n}\begin{pmatrix} \hat{\theta}_n - \theta_0 \\ t_{\hat{\theta}_n} - 0 \end{pmatrix} \xrightarrow{d} \mathcal{N}\left(0, \begin{bmatrix} V_{\theta_0} & 0 \\ 0 & U_{\theta_0} \end{bmatrix}\right),$$

*where $V_{\theta_0} = 4\pi(D'_{\theta_0} W^{-1}_{\theta_0} D_{\theta_0})^{-1}$ and $U_{\theta_0} = \pi^{-1} W^{-1}_{\theta_0}(I_{r \times r} - (4\pi)^{-1} D_{\theta_0} V_{\theta_0} \times D'_{\theta_0} W^{-1}_{\theta_0})$.*

(ii) *Additionally, suppose* A.5 *and $f = f_{\theta_0}$ hold, that $\partial f_{\theta_0}/\partial \theta$ is Riemann integrable, that $P(0 \in \operatorname{ch}^{\circ}\{\pi G_\theta(\lambda_j)[I_n(\lambda_j) - f_\theta(\lambda_j)]\}_{j=1}^N, \theta \in \overline{\mathcal{B}}_n) \to 1$ for the closure $\overline{\mathcal{B}}_n$, and that in a neighborhood of $\theta_0$, $\int_0^\pi G_\theta(\lambda) f_\theta(\lambda)\, d\lambda = \mathcal{M}$ and $\|\partial f_\theta/\partial \theta\|, \|\partial^2 f_\theta(\lambda)/\partial\theta\partial\theta'\| \leq C|\lambda|^{-\alpha}$, $\lambda \in (0, \pi]$. Then there exists a sequence of statistics $\{\hat{\theta}_{n,\mathcal{F}}\}$ such that $P$ ($\hat{\theta}_{n,\mathcal{F}}$ is a maximum of $R_{n,\mathcal{F}}(\theta)$ and $\hat{\theta}_{n,\mathcal{F}} \in \mathcal{B}_n$) $\to 1$ as $n \to \infty$ and the distributional result in* (i) *is valid for $\sqrt{n}(\hat{\theta}_{n,\mathcal{F}} - \theta_0, t_{\hat{\theta}_{n,\mathcal{F}}}/2)'$.*

(iii) *If $\kappa_{4,\varepsilon} = 0$, then Theorem 2(ii) holds, even if $\mathcal{M} \neq 0 \in \mathbb{R}^r$.*

REMARK 2. When assuming $f \in \mathcal{F}$ in Theorem 2(ii), a constant function $\int_0^\pi G_\theta(\lambda) f_\theta(\lambda)\, d\lambda = \mathcal{M}$ of $\theta$ represents a natural relationship between the chosen estimating functions and $\mathcal{F}$ [e.g., the Whittle estimating equations $G_\theta^w$ in (11)]. The probabilistic assumption on the closure $\overline{\mathcal{B}}_n$ in Theorem 2(ii) is similar to Assumption (A.5) and implies the FDEL ratio $\ell_{n,\mathcal{F}}(\theta)$ exists finitely in a neighborhood of $\theta_0$.

As pointed out by a referee, Theorem 2 establishes consistency of a local maximizer of the EL function only. In the event that the likelihood $R_n$ or $R_{n,\mathcal{F}}$ has a single maximum with probability approaching 1 [e.g., by concavity of $R_n(\theta)$], the sequence $\{\hat{\theta}_n\}$ or $\{\hat{\theta}_{n,\mathcal{F}}\}$ corresponds to a global MELE. In many cases, the consistency of global maximizers can also be established using additional conditions. In Theorem 3 below, we give conditions for the consistency of $\hat{\theta}_n$; similar conditions for $\hat{\theta}_{n,\mathcal{F}}$ can be developed. Note that these conditions are satisfied by the estimating functions given in Section 2.

THEOREM 3. *Suppose Assumption* A.1 *holds.*

(i) *Assume $\theta_0$ lies in the interior of $\Theta$, $G_\theta(\lambda)$ is a (componentwise) continuous and monotone function of $\theta$ for $\lambda \in \Pi$ and, for $\theta \in \Theta$, $|G_\theta(\lambda)|$ is Riemann integrable and bounded by $C|\lambda|^{-\delta}$ for some $\delta < 1$, $\alpha + \delta < 1$. Then as $n \to \infty$,*

(18) $\qquad P\left(\hat{\theta}_n = \arg\max_{\theta \in \Theta} R_n(\theta) \text{ exists}\right) \longrightarrow 1 \quad \text{and} \quad \hat{\theta}_n \xrightarrow{p} \theta_0.$

(ii) *Suppose $\Theta$ is compact, $G_\theta(\lambda) \equiv G_\theta^w$ from* (11) *[or $G_\theta(\lambda) \equiv G_\vartheta^{w*}(\lambda)$ from* (12)*] is continuous at all $(\lambda, \theta) \in \Pi \times \Theta$ and $W(\theta)$ from* (8) *attains its minimum on the interior of $\Theta$. Then* (18) *holds as $n \to \infty$.*

5.2. *Empirical likelihood tests of hypotheses.* EL ratio test statistics with $\hat{\theta}_n$ and $\hat{\theta}_{n,\mathcal{F}}$ are possible for both parameter and moment hypotheses. Similarly to parametric likelihood, we can use the log-EL ratio $\ell_n(\theta_0) - \ell_n(\hat{\theta}_n)$ to test the parameter assumption $H_0: \theta = \theta_0$. For testing the null hypothesis that the true parameter $\theta_0$ satisfies the spectral mean condition in (6), the log-ratio statistic $\ell_n(\hat{\theta}_n)$ is useful. Analogous tests are possible with $\ell_{n,\mathcal{F}}(\theta_0)$ and $\ell_{n,\mathcal{F}}(\hat{\theta}_{n,\mathcal{F}})$. These EL log-ratio statistics have limiting chi-square distributions for testing the above hypotheses.

THEOREM 4. *Under the assumptions of Theorem* 2 *with the sequences* $\{\hat{\theta}_n\}$ *and* $\{\hat{\theta}_{n,\mathcal{F}}\}$,

(i) $\ell_n(\theta_0) - \ell_n(\hat{\theta}_n) \xrightarrow{d} \chi_p^2$, $\ell_n(\hat{\theta}_n) \xrightarrow{d} \chi_{r-p}^2$ *and these are asymptotically independent*;

(ii) $\ell_{n,\mathcal{F}}(\theta_0) - \ell_{n,\mathcal{F}}(\hat{\theta}_{n,\mathcal{F}}) \xrightarrow{d} \chi_p^2$, $\ell_{n,\mathcal{F}}(\hat{\theta}_{n,\mathcal{F}}) \xrightarrow{d} \chi_{r-p}^2$ *and these are asymptotically independent, if the assumptions in Theorem* 2(ii) *are additionally satisfied*;

(iii) *if* $\kappa_{4,\varepsilon} = 0$, *Theorem* 4(ii) *remains valid, even if* $\mathcal{M} \neq 0 \in \mathbb{R}^r$.

In Sections 6 and 7, Theorems 1–4 will be applied to Whittle estimation and goodness-of-fit testing in the FDEL framework.

5.2.1. *Parameter restrictions and nuisance parameters.* Qin and Lawless [39] introduced constrained EL inference for independent samples and [23] provided a blockwise version for time domain EL under SRD. We can also consider FDEL estimation subject to the following system of parameter constraints on $\theta: \psi(\theta) = 0 \in \mathbb{R}^q$, where $q < p$ and $\Psi(\theta) = \partial \psi(\theta)/\partial \theta$ is of full row rank $q$. By maximizing the EL functions in (15) or (16) under the above restrictions, we find constrained MELE's $\hat{\theta}_n^\psi$ or $\hat{\theta}_{n,\mathcal{F}}^\psi$.

COROLLARY 1. *Suppose the conditions in Theorem* 2 *hold and, in a neighborhood of* $\theta_0$, $\psi(\theta)$ *is continuously differentiable,* $\|\partial^2 \psi(\theta)/\partial \theta \partial \theta'\|$ *is bounded and* $\Psi(\theta_0)$ *is rank* $q$. *If* $H_0: \psi(\theta_0) = 0$ *holds, then* $\ell_n(\hat{\theta}_n^\psi) - \ell_n(\hat{\theta}_n) \xrightarrow{d} \chi_q^2$ *and* $\ell_n(\theta_0) - \ell_n(\hat{\theta}_n^\psi) \xrightarrow{d} \chi_{p-q}^2$ *as* $n \to \infty$.

We can then sequentially test $H_0: \psi(\theta_0) = 0$ with a log-likelihood ratio statistic $\ell_n(\hat{\theta}_n^\psi) - \ell_n(\hat{\theta}_n)$ and, if failing to reject $H_0$, make an approximate $100(1-\gamma)\%$ confidence region for constrained $\theta$ values $\{\theta: \psi(\theta) = 0, \ell_n(\theta) - \ell_n(\hat{\theta}_n^\psi) \leq \chi_{p-q,1-\gamma}^2\}$.

Profile FDEL ratio statistics can also be developed to conduct tests in the presence of nuisance parameters (see Corollary 5 of [38] for the i.i.d. case).

Suppose $\theta = (\theta_1', \theta_2')'$, where $\theta_1$ and $\theta_2$ are $q \times 1$ and $(p-q) \times 1$ vectors, respectively. For fixed $\theta_1 = \theta_1^0$, suppose that $\hat{\theta}_2^0$ and $\hat{\theta}_{2,\mathcal{F}}^0$ maximize the EL functions $R_n(\theta_1^0, \theta_2)$ and $R_{n,\mathcal{F}}(\theta_1^0, \theta_2)$ with respect to $\theta_2$.

COROLLARY 2. *Under the conditions in Theorem 2, if* $H_0 : \theta_1 = \theta_1^0$ *holds, then* $\ell_n(\theta_1^0, \hat{\theta}_2^0) - \ell_n(\hat{\theta}_n) \xrightarrow{d} \chi_q^2$ *as* $n \to \infty$.

If the assumptions in Theorem 2(ii) are also satisfied, Corollaries 1 and 2 hold using $\ell_{n,\mathcal{F}}(\cdot)$, $\hat{\theta}_{n,\mathcal{F}}$, $\hat{\theta}_{n,\mathcal{F}}^\psi$ and $\hat{\theta}_{2,\mathcal{F}}^0$.

## 6. Whittle estimation.

EXAMPLE 4 (Continued). With SRD linear processes, Monti [31] suggested EL confidence regions for Whittle-like estimation of parameters $\theta = (\sigma^2, \vartheta')'$ characterizing $f_\theta \in \mathcal{F}$ from (9). Theorem 1 provides two refinements of [31].

*Refinement 1.* Monti [31] develops an EL function for $\theta$ by treating the standardized ordinates $I_n(\lambda_j)/f_\theta(\lambda_j)$, $j = 1, \ldots, N$, as approximately i.i.d. random variables, similarly to the FDB. The EL ratio in (4.1) of [31] asymptotically corresponds to $\ell_{n,\mathcal{F}}(\theta)$ when using the estimating functions $G_\theta^w$ from (11). With this choice of functions, the nonzero spectral mean $\mathcal{M}_w \neq 0$ in (11) is due to the first estimating function $f_\theta^{-1}$ intended to prescribe $\sigma^2$. Note that the use of $\ell_{n,\mathcal{F}}(\theta)$ and $G_\theta^w$ for setting confidence regions *requires* the additional assumption that the fourth order innovation cumulant $\kappa_{4,\varepsilon} = 0$, by Theorem 1. Valid joint confidence regions are otherwise not possible here because $\mathcal{M}_w \neq 0$. This complication due to $\sigma^2$-inference is related to the inconsistency of the FDB Whittle estimate of $\sigma^2$ when $\kappa_{4,\varepsilon} \neq 0$, as described by [11]. While valid for SRD Gaussian series with $\kappa_{4,\varepsilon} = 0$, periodogram-based EL formulation in [31] may not be applicable to general SRD linear processes.

*Refinement 2.* Treating $\sigma^2$ as a nuisance parameter and concentrating it out of the Whittle likelihood, Monti [31] suggests an EL ratio statistic for estimation of the remaining $p-1$ parameters $\vartheta$ in (9) via confidence regions. The statistic (6.1) of [31] appears to be asymptotically equivalent to $-2\log R_n(\vartheta) = 1/2 \cdot \ell_n(\vartheta)$, based on the $p-1$ estimating functions $G_\vartheta^{w*}$ from (12) and the $N = \lfloor (n-1)/2 \rfloor$ periodogram ordinates. Note that for the function $G_\vartheta^{w*}$, we have $\mathcal{M}_{w*} = 0$ so that Theorem 1(i) applies to $\ell_n(\vartheta)$.

For Whittle-like estimation of $\vartheta$ in the parameterization from (9), the EL log-ratio $\ell_n(\vartheta)$ based on the functions $G_\vartheta^{w*}$ in (12) appears preferable to $\ell_{n,\mathcal{F}}(\vartheta)$. This selection results in asymptotically correct confidence regions for $\vartheta$ under both SRD and LRD, even for misspecified situations ($f \notin \mathcal{F}$) where the moments in (10) still hold.

## 7. Goodness-of-fit tests.

### 7.1. *Simple hypothesis case.*

EXAMPLE 3 (Continued). We return to the simple hypothesis test $H_0: f = f_0$ for some possible density $f_0$. To assess the goodness-of-fit, Milhoj [30] and Beran [3] proposed the test statistic $T_n = \pi A_n / B_n^2$ for mixing SRD linear processes (with $\kappa_{4,\varepsilon} = 0$) and long-memory Gaussian processes, respectively, and established the limiting bivariate normal law of $\sqrt{n}\{(A_n, B_n)' - (2\pi, \pi)'\}$ under $H_0$ for $A_n = 2\pi/n \sum_{j=1}^{N} I_n^2(\lambda_j)/f_0^2(\lambda_j)$, $B_n = 2\pi/n \sum_{j=1}^{N} I_n(\lambda_j)/f_0(\lambda_j)$. Since these linear processes involve $\kappa_{4,\varepsilon} = 0$, under Theorem 1(iii), we can construct a single statistic $\ell_{n,\mathcal{F}}$ to test $H_0$ by treating $\mathcal{F} = \{f_0\}$ in (7) and employing a single estimating function $f_0^{-1}$ satisfying (6) with $\mathcal{M} = \pi$ under $H_0$. In expanding $\ell_{n,\mathcal{F}} = n(B_n - \pi)^2/(\pi A_n) + o_p(1)$, we find the FDEL ratio statistic asymptotically incorporates much of the same information in $T_n$ under $H_0$ (with better power when $f = cf_0$, $c \neq 1$).

For testing the special hypothesis $H_0: \{X_t\}$ is white noise (constant $f$), a FDEL goodness-of-fit test based on the process autocorrelations can be applied, similar to some Portmanteau tests [27, 28]. The estimating functions $G(\lambda) = (\cos(\lambda), \ldots, \cos(m\lambda))'$ satisfy (6) with $\mathcal{M} = 0 \in \mathbb{R}^m$ under this $H_0$ (see Example 1) and yield a single EL ratio $\ell_n$ which pools information across $m$ EL-estimated autocorrelation lags.

### 7.2. *Composite hypothesis case.* 

To test the composite hypothesis $H_0: f \in \mathcal{F}$ for a specified parametric class $\mathcal{F}$, various frequency domain tests have been proposed by [3, 30, 37] which use Whittle estimation to select the "best" fitting model in $\mathcal{F}$ and then compare the fitted density to the periodogram across all ordinates. We show that FDEL techniques can produce similar goodness-of-fit tests, while expanding our EL theory slightly.

Suppose $\{X_t\}$ is a Gaussian time series and that we wish to test if $f \in \mathcal{F}$ for some parametric family as in (9), which includes densities (3) or (4). This scenario is considered by [3] and [37] for LRD and SRD Gaussian models, respectively. FDEL methods may simultaneously incorporate both components of model fitting and model comparison through estimating equations

$$(19) \qquad \int_0^\pi G_\theta^w f \, d\lambda = \mathcal{M}_w, \qquad \int_0^\pi (f/f_\theta)^2 \, d\lambda = \pi,$$

where $G_\theta^w = (f_\theta^{-1}, \partial f_\theta^{-1}/\partial \vartheta')'$ are the Whittle estimating functions from (11) for the parameters $\theta = (\sigma^2, \vartheta')' \in \mathbb{R}^p$ in $f_\theta$. Note that we introduce an overidentifying moment restriction on $f$ in (19) so that $r = p + 1$. We then extend the log-likelihood statistic $\ell_{n,\mathcal{F}}$ in (17) to include $I_n^2$ ordinates by defining $\ell_{I_n^2,\mathcal{F}}(\theta) = -2\log R_{I_n^2,\mathcal{F}}(\theta)$ for $R_{I_n^2,\mathcal{F}}(\theta) = (N/\pi)^N L_{I_n^2,\mathcal{F}}(\theta)$ and

$$L_{I_n^2,\mathcal{F}}(\theta) = \sup\Bigg\{\prod_{j=1}^N w_j : w_j \geq 0, \ \sum_{j=1}^N w_j = \pi,$$
$$\sum_{j=1}^N w_j \begin{pmatrix} I_n^2(\lambda_j)/\{2f_\theta^2(\lambda_j)\} - 1 \\ G_\theta^w(\lambda_j)[I_n(\lambda_j) - f_\theta(\lambda_j)] \end{pmatrix} = 0\Bigg\},$$

using $f_\theta$, $f_\theta^2$ above to approximate $\mathrm{E}(I_n)$, $\mathrm{E}(I_n^2/2)$ for each ordinate.

To evaluate $H_0: f \in \mathcal{F}$, we test if the moment conditions in (19) hold for some $\theta$ value. Following the test prescribed by Theorem 4, we find the argument maximum of $L_{I_n^2,\mathcal{F}}(\theta)$, say $\hat{\theta}_{I_n^2,\mathcal{F}}$, and form a test statistic $\ell_{I_n^2,\mathcal{F}}(\hat{\theta}_{I_n^2,\mathcal{F}})$ for this $H_0$. The subsequent extension of Theorem 4 gives the distribution for our test statistic.

THEOREM 5. *Suppose $\{X_t\}$ is Gaussian and the assumptions in Theorem 4(ii) hold for $f_{\theta_0}$ and $G_{\theta_0}^w$ with $\alpha - \beta < \upsilon$ for each arbitrarily small $\upsilon > 0$. Under the null hypothesis $f \equiv f_{\theta_0} \in \mathcal{F}$, $\ell_{I_n^2,\mathcal{F}}(\hat{\theta}_{I_n^2,\mathcal{F}}) \xrightarrow{d} \chi_1^2$ as $n \to \infty$.*

The distributional result is valid even with nonzero spectral mean conditions in (19) because the process is Gaussian. We make a few comments about model misspecification. Suppose $f \notin \mathcal{F}$, but $\theta_0$ still represents the parameter value which minimizes the asymptotic distance measure $W(\theta)$ in (8) and $f_{\theta_0}$ satisfies (10) (i.e., $\int_0^\pi G_{\theta_0}^w f\,d\lambda = \mathcal{M}_w$ holds). A test consistency property can then be established: as $n \to \infty$,

$$n^{-1}\ell_{I_n^2,\mathcal{F}}(\hat{\theta}_{I_n^2,\mathcal{F}}) \xrightarrow{p} a_0\bigg\{\int_0^\pi \bigg(\frac{f(\lambda)}{f_{\theta_0}(\lambda)} - 1\bigg)^2 d\lambda\bigg\}^2 > 0,$$

where $a_0 > 0$ depends on $f$ and $f_{\theta_0}$. We are assured that the test statistic can determine if $H_0: f \in \mathcal{F}$ is true as the sample size increases.

**8. Conclusions.** We have introduced a frequency domain version of empirical likelihood (FDEL) based on the periodogram, which allows spectral inference, in a variety of applications, on both short- and long-range dependent linear processes. Further numerical study and development of FDEL will be considered in future communications. See [32] for extensions of FDEL using the tapered periodogram in (13) under weak dependence. A valuable area of potential research includes Bartlett corrections to the EL ratios in (17) (see [12] for i.i.d. data). Second- and higher-order correct FDEL confidence regions may be possible without the kernel estimation or stringent moment assumptions required with the frequency domain bootstrap of [11].

**9. Proofs.** We only outline some proofs here for reasons of brevity. Detailed proofs can be found in [32] and we shall refer to relevant results given there. We require some additional notation and lemmas to facilitate the proofs. In the following, $C$ or $C(\cdot)$ will denote generic constants that depend on their arguments (if any), but do not depend on $n$, including ordinates $\{\lambda_j\}_{j=1}^N$.

Define the mean-corrected discrete Fourier transforms $d_{nc}(\lambda) = \sum_{t=1}^n (X_t - \mu)e^{-\imath t\lambda}$, $\lambda \in \Pi$. Note that $2\pi n I_{nc}(\lambda) = |d_{nc}(\lambda)|^2 = d_{nc}(\lambda)d_{nc}(-\lambda)$ and $I_{nc}(\lambda_j) = I_n(\lambda_j)$ for $j = 1, \ldots, N$. Let $H_n(\lambda) = \sum_{t=1}^n e^{-\imath t\lambda}$, $\lambda \in \mathbb{R}$, and write $K_n(\lambda) = (2\pi n)^{-1}|H_n(\lambda)|^2$ to denote the Fejér kernel. The function $K_n$ is nonnegative and even with period $2\pi$ on $\mathbb{R}$ and satisfies $\int_\Pi K_n \, d\lambda = 1$ (see page 71 of [6]). We adopt the standard that an even function $g : \Pi \longrightarrow \mathbb{R}$ can be periodically (period $2\pi$) extended to $\mathbb{R}$ by $g(\lambda) = g(-\lambda)$, $g(\lambda) = g(\lambda + 2\pi)$ for $\lambda \in \mathbb{R}$. We make extensive use of the following function from [8]. Let $L_{ns} : \mathbb{R} \longrightarrow \mathbb{R}$ be the periodic extension of

$$L_{ns}(\lambda) =: \begin{cases} e^{-s}n, & |\lambda| \leq e^s/n, \\ \dfrac{\log^s(n|\lambda|)}{|\lambda|}, & e^s/n < |\lambda| \leq \pi, \end{cases} \quad \lambda \in \Pi, \ s = 0, 1.$$

For each $n \geq 1$, $s = 0, 1$, $L_{ns}(\cdot)$ is decreasing on $[0, \pi]$ and

(20) $\quad |H_n(\lambda)| \leq C L_{n0}(\lambda), \quad L_{n0}(\lambda) \leq \dfrac{3\pi n}{1 + |\lambda \bmod 2\pi|n}, \quad \lambda \in \mathbb{R}.$

In the following, $\operatorname{cum}(Y_1, \ldots, Y_m)$ denotes the joint cumulant of random variables $Y_1, \ldots, Y_m$. We often refer to cumulant properties from Section 2.3 of [5], including the product theorem for cumulants (cf. Theorem 2.3.2).

We remark that Lemma 5 to follow ensures that the log-likelihood ratio $\ell_n(\theta_0)$ exists asymptotically. Lemmas 6 and 7 establish important results for Riemann integrals based on the periodogram under both LRD and SRD; Lemma 6 considers the distribution of empirical spectral means and Lemma 7 is required for variance estimation.

LEMMA 1. *Let $1 \leq i \leq j \leq N$ and $0 < d < 1$. If $a_i \in \{\pm \lambda_i\}$, $a_j \in \{\pm \lambda_j\}$ and $a_i + a_j \neq 0$, then*

(i) $\quad L_{n0}(a_i + a_j) \leq \dfrac{nc_{ijn}}{2\pi},$
$$c_{ijn} = \begin{cases} (j-i)^{-1}, & \operatorname{sign}(a_i) \neq \operatorname{sign}(a_j), \\ (j+i)^{-1}, & \operatorname{sign}(a_i) = \operatorname{sign}(a_j), \ i+j \leq n/2, \\ (n-j-i)^{-1}, & \operatorname{sign}(a_i) = \operatorname{sign}(a_j), \ i+j > n/2. \end{cases}$$

(ii) $\quad L_{n1}(a_i + a_j) \leq \min\{\log(n\pi)L_{n0}(a_i + a_j), \ nC(d)[c_{ijn}]^d\}.$

(iii) $$\int_\Pi L_{n0}^m(\lambda)\,d\lambda \leq C(m)[\log(n) + n^{m-1}], \qquad m \geq 1 \in \mathbb{Z}.$$

(iv) $$\int_\Pi L_{n0}^m(r_1+\lambda)L_{n0}^m(r_2-\lambda)\,d\lambda \leq \begin{cases} CL_{n1}(r_1+r_2), & m=1, \\ CnL_{n0}^2(r_1+r_2), & m=2, \end{cases} \qquad r_1, r_2 \in \mathbb{R}.$$

PROOF. Parts (iii)–(iv) are from Lemmas 1 and 2 of [8]. Lemma 1(i) follows from the fact that $|(a_i + a_j) \bmod 2\pi| \geq 2\pi/n$ if $a_i + a_j \neq 0$, along with the definition of $L_{n0}$. □

LEMMA 2. *Suppose Assumption* A.1 *holds. Let* $\Pi_\rho = [\rho, \pi]$ *for* $0 < \rho < \pi$. *If* $a_1, a_2 \in \Pi$, $|a_1| \leq |a_2|$ *and* $|a_2| \in \Pi_\rho$, *then* $\mathrm{cum}(d_{nc}(a_1), d_{nc}(a_2)) = 2\pi H_n(a_1+a_2)f(a_2) + R_{n\rho}(a_1, a_2)$ *and* $R_{n\rho} = o(n)$ *holds for* $R_{n\rho} \equiv \sup\{|R_{n\rho}(a_1, a_2)| : a_1, a_2 \in \Pi, |a_1| \leq |a_2|, |a_2| \in \Pi_\rho\}$.

PROOF. We modify the proof of Theorem 1(a) of [8]; see [32], Lemma 5 for details. □

LEMMA 3. *Let* $1 \leq i \leq j \leq N$ ($n \geq 3$) *and* $a_1, \ldots, a_k \in \{\pm\lambda_i, \pm\lambda_j\}$, $|a_1| \leq |a_2| \leq \cdots \leq |a_k|$, *with* $2 \leq k \leq 8$. *Under Assumption* A.1,

(i) $$|\mathrm{cum}(d_{nc}(a_1), d_{nc}(a_2))| \leq C \begin{cases} |a_1|^{-\alpha}(|a_2|^{-1} + L_{n1}(a_1+a_2)), & \text{if } \alpha > 0, \\ L_{n1}(a_1+a_2), & \text{if } \alpha = 0. \end{cases}$$

(ii) $$|\mathrm{cum}(d_{nc}(a_1), \ldots, d_{nc}(a_k))| \leq C\{|a_k|^{\alpha/2-1}|a_{k-1}|^{-1/2} + n\log^{k-1}(n)\} \prod_{j=1}^k |a_j|^{-\alpha/2}.$$

PROOF. We show Lemma 3(ii); the proof of (i) can be similarly shown with details given in [32], Lemma 1. From [25, 41], the $k$th order joint cumulant ($2 \leq k \leq 8$) may be expressed as $\mathrm{cum}(d_{nc}(a_1), \ldots, d_{nc}(a_k)) = (2\pi)^{-k+1}\kappa_{\varepsilon,k}\nu(\Pi^{k-1})$ using the $k$th order innovation cumulant $\kappa_{\varepsilon,k}$ and a function $\nu$ of Borel measurable sets $A \subset \Pi^{k-1}$ defined as

$$\nu(A) = \int_A H_n\left(\sum_{j=1}^k a_j - \sum_{j=1}^{k-1} z_j\right) b\left(\sum_{j=1}^{k-1} a_j - z_j\right) \\ \times \prod_{j=1}^{k-1} \{H_n(z_j)b(z_j - a_j)\}\,dz_1 \ldots dz_{k-1}.$$

On $B = \bigcap_{j=1}^{k-1} \{(z_1,\ldots,z_{k-1}) \in \Pi^{k-1} : |z_j - a_j| \leq |a_j|/2k\}$, we have $|H_n(z_j)| \leq C|a_j|^{-1}$, $|H_n(\sum_{j=1}^{k} a_j - \sum_{j=1}^{k-1} z_j)| \leq C|a_k|^{-1}$ by (20) and, by applying Holder's inequality,

$$\int_{\substack{|z_{k-1}-a_{k-1}| \\ \leq |a_{k-1}|/2k}} \left| b\left(\sum_{j=1}^{k-1} a_j - z_j\right) b(z_{k-1} - a_{k-1}) \right| dz_{k-1} \leq C|a_{k-1}|^{(1-\alpha)/2},$$

while $\int_{|z_j - a_j| \leq |a_j|/2k} |b(z_j - a_j)| dz_j \leq C|a_j|^{1-\alpha/2}$ for $1 \leq j \leq k-2$. By these inequalities, $|\nu(B)| \leq C|a_k|^{-1}|a_{k-1}|^{-1/2} \prod_{j=1}^{k-1} |a_j|^{-\alpha/2}$. Now, on $B_j = \{(z_1,\ldots,z_{k-1}) \in \Pi^{k-1} : |z_j - a_j| > |a_j|/2k\}$ for fixed $1 \leq j \leq k-1$, we find $|b(z_j - a_j)| \leq C|a_j|^{-\alpha/2}$, $|H_n(z_j)| \leq Cn$ and

$$\int_\Pi \left| H_n\left(\sum_{\ell=1}^k a_\ell - \sum_{\ell=1}^{k-1} z_\ell\right) b\left(\sum_{\ell=1}^{k-1} a_\ell - z_\ell\right) \right| dz_j$$
$$\leq \int_{|\lambda| \leq |a_k|/2} C|a_k|^{-1} |\lambda|^{-\alpha/2} d\lambda + \int_{|\lambda| > |a_k|/2} C|a_k|^{-\alpha/2} L_{n0}(a_k - \lambda) d\lambda$$
$$\leq C|a_k|^{-\alpha/2} \log(n),$$

by Lemma 1(iii), while $\int_\Pi |H_n(z_\ell)||b(z_\ell - a_\ell)| dz_\ell \leq C|a_\ell|^{-\alpha/2} \log(n)$ for each $\ell \neq j$. Since $\Pi^{k-1} \setminus B = \bigcup_{j=1}^{k-1} B_j$, we have $|\nu(\Pi^{k-1} \setminus B)| \leq \sum_{j=1}^{k-1} |\nu(B_j)| \leq Cn \log^{k-1}(n) \prod_{j=1}^{k} |a_j|^{-\alpha/2}$. $\square$

LEMMA 4. *Under Assumption* A.1, $\hat{r}(k) = \hat{r}(-k) = n^{-1} \sum_{j=1}^{n-k} (X_j - \mu) \times (X_{j+k} - \mu) \xrightarrow{p} r(k) = \text{Cov}(X_j, X_{j+k})$ *as* $n \to \infty$ *for each* $k \geq 0$.

PROOF. $\hat{r}(k)$ is asymptotically unbiased and $\text{Var}(\hat{r}(k)) = o(1)$, by Lemma 3.3 of [22]. $\square$

LEMMA 5. *Under Assumption* A.1, *suppose* $G = (g_1,\ldots,g_r)' \equiv G_{\theta_0}$ *is even with finite discontinuities on* $[0,\pi]$ *and satisfies Assumptions* A.2 *and* A.4. *If* $\int_\Pi Gf \, d\lambda = 0 \in \mathbb{R}^r$, *then* $P(0 \in \text{ch}°\{\pi G(\lambda_j) I_n(\lambda_j)\}_{j=1}^N) \to 1$ *as* $n \to \infty$.

PROOF. See [32] for details. It can be shown that $\inf_{\|y\|=1} \int_\Pi \mathbb{1}_{\{y'G(\lambda)>0\}} y' \times Gf \, d\lambda \geq C > 0$ and that $P(\inf_{\|y\|=1} \frac{4\pi}{n} \sum_{j=1}^{N} I_n(\lambda_j) y' G(\lambda_j) \mathbb{1}_{\{y'G(\lambda_j)>0\}} \geq C/2) \to 1$ follows using Lemma 4 with arguments from Lemma 1 of [20]. When the event in this probability statement holds, the separating/supporting hyperplane theorem implies that $0 \in \text{ch}°\{\pi G(\lambda_j) I_n(\lambda_j)\}_{j=1}^N$. $\square$

LEMMA 6. *Suppose Assumptions* A.1–A.3 *hold with respect to an even function* $G = (g_1, \ldots, g_r)' \equiv G_{\theta_0}$ *and let* $J_n = (2\pi/n) \sum_{j=1}^{N} G(\lambda_j) I_n(\lambda_j)$. *Then* $\sqrt{n}(J_n - \int_0^\pi fG\, d\lambda) \xrightarrow{d} \mathcal{N}(0, V)$ *as* $n \to \infty$, *where*

$$V = \pi \int_\Pi f^2 GG' \, d\lambda + \frac{\kappa_{4,\varepsilon}}{4\sigma_\varepsilon^4} \left( \int_\Pi fG\, d\lambda \right) \left( \int_\Pi fG\, d\lambda \right)'.$$

*If, in addition,* A.5 *holds, then* $\sqrt{n}\tilde{J}_n \xrightarrow{d} \mathcal{N}(0, V)$ *for* $\tilde{J}_n = (2\pi/n) \sum_{j=1}^{N} G(\lambda_j) \times [I_n(\lambda_j) - f(\lambda_j)]$.

PROOF. By Assumptions A.1–A.2, we have that $\sqrt{n}(\int_0^\pi GI_{nc}\, d\lambda - \mathrm{E} \int_0^\pi GI_{nc}\, d\lambda) \xrightarrow{d} N(0, V)$ from Theorem 2 and Lemma 6 of [17] and the Cramer–Wold device. To show the result for $J_n$ in Lemma 6, it suffices to establish that

(21)
$$\left\| \mathrm{E} \int_0^\pi G(I_{nc} - f)\, d\lambda \right\| = o(n^{-1/2}),$$
$$\left\| J_n - \int_0^\pi GI_{nc}\, d\lambda \right\| = o_p(n^{-1/2}).$$

The distribution of $\tilde{J}_n$ then also follows, using $\|(2\pi/n) \sum_{j=1}^{N} G(\lambda_j) f(\lambda_j) - \int_0^\pi fG\, d\lambda\| = o(n^{-1/2})$, which can be established with straightforward arguments; see [32], Lemma 10.

Without loss of generality, we assume that $G = g$ (i.e., $r = 1$) and establish (21) under Condition 2 of Assumption A.3; see [32] for proofs under Conditions 1 or 3 of A.3. Proving (21) under Condition 1 of A.3 involves using the $n$th Cesàro mean $c_n g(\lambda) = \int_\Pi K_n(\lambda - y) g(y)\, dy$, $\lambda \in \Pi$, for which $\sup_{\lambda \in \Pi} |g(\lambda) - c_n g(\lambda)| = o(n^{-1/2})$, by Theorem 6.5.3 of [14], and using $\frac{2\pi}{n} \sum_{j=-N}^{\lfloor n/2 \rfloor} c_n g(\lambda_j) I_{nc}(\lambda_j) = \int_\Pi c_n GI_{nc}\, d\lambda$, $\mathrm{E} \int_\Pi gI_{nc}\, d\lambda = \int_\Pi c_n gf\, d\lambda$. The proof under Condition 3 relies on Theorem 3.2 and Lemma 3.1 of [9] and Lemma 4 of [8].

We show the first convergence in (21) here; Lemma A.1 in the Appendix gives the second part of (21). Using the evenness of $K_n$, $\int_\Pi K_n(\lambda - y)\, dy = 1$ and $\mathrm{E}(I_{nc}(\lambda)) = \int_\Pi K_n(\lambda - y) f(y)\, dy$, we have

(22)
$$\sqrt{n} \left| \mathrm{E} \int_\Pi gI_{nc}\, d\lambda - \int_\Pi gf\, d\lambda \right| \leq \sqrt{n} \int_{\Pi^2} K_n(\lambda - y) f(y) |g(\lambda) - g(y)|\, dy\, d\lambda$$
$$\leq C n^{3/2} \int_{(0,\pi]^2} \frac{f(y)|g(\lambda) - g(y)|}{(1 + |\lambda - y|n)^2}\, dy\, d\lambda,$$

where the last inequality follows from (20) and $|(\lambda - y) \mod 2\pi| \geq ||\lambda| - |y||$, $\lambda, y \in \Pi$. We now modify an argument from [17] (page 99). Under Condition 2 of A.3 (with respect to $\beta_1 > 0$), we may pick $0 < \gamma < 1/2$ so that $0 < \gamma^* \equiv \gamma + \beta_1(1-\gamma) - \alpha < 1$ and $f(y)|g(\lambda) - g(y)| \leq C(\min\{y, \lambda\})^{-1+\gamma^*}|\lambda - y|^{1-\gamma}$, $\lambda, y \in (0, \pi]$. We then bound (22) by

$$Cn^{-1/2-\gamma^*+\gamma} \int_0^{n\pi} \left( \int_0^\infty \frac{y^{-1+\gamma^*}}{(1+|\lambda-y|)^{1+\gamma}} \, dy \right) d\lambda \leq Cn^{-1/2+\gamma} = o(1),$$

using the fact that there exist a $C > 0$ such that $\int_0^\infty y^{-1+\gamma^*}(1+|\lambda-y|)^{-1-\gamma} \, dy \leq C|\lambda|^{-1+\gamma^*}$ for any $\lambda \in \mathbb{R}$. $\square$

LEMMA 7. *Under Assumption* A.1, *suppose* $g$ *and* $h$ *are real-valued, even Riemann integrable functions on* $\Pi$ *such that* $|g(\lambda)|, |h(\lambda)| \leq C|\lambda|^\beta$, $0 \leq \beta < 1$, $\alpha - \beta < 1/2$. *Then as* $n \to \infty$,

$$\frac{2\pi}{n} \sum_{j=1}^N g(\lambda_j) h(\lambda_j) I_n^2(\lambda_j) \text{ and } 2 \cdot \frac{2\pi}{n} \sum_{j=1}^N g(\lambda_j) h(\lambda_j) (I_n(\lambda_j) - f(\lambda_j))^2$$
$$\xrightarrow{p} \int_\Pi ghf^2 \, d\lambda.$$

PROOF. We consider the first Riemann sum above; convergence of the second sum can be similarly established. Since $(2\pi/n) \sum_{j=1}^N g(\lambda_j) h(\lambda_j) f^2(\lambda_j) \to \int_0^\pi ghf^2 \, d\lambda$ by the Lebesgue dominated convergence theorem, it suffices to establish

$$B_n \equiv \left| \mathrm{E}(S_n) - \frac{4\pi}{n} \sum_{j=1}^N g(\lambda_j) h(\lambda_j) f^2(\lambda_j) \right| = o(1), \qquad \mathrm{Var}(S_n) = o(1),$$

for $S_n = (2\pi/n) \sum_{j=1}^N g(\lambda_j) h(\lambda_j) I_n^2(\lambda_j)$. We show that $B_n = o(1)$; Lemma A.2 in the Appendix shows that $\mathrm{Var}(S_n) = o(1)$. By $\mathrm{E}(d_{nc}(\lambda)) = 0$ and the product theorem for cumulants,

$$(2\pi n)^2 \mathrm{E}(I_n^2(\lambda_j)) = \mathrm{cum}^2(d_{nc}(\lambda_j), d_{nc}(\lambda_j)) + 2\mathrm{cum}^2(d_{nc}(\lambda_j), d_{nc}(-\lambda_j))$$
$$+ \mathrm{cum}(d_{nc}(\lambda_j), d_{nc}(\lambda_j), d_{nc}(-\lambda_j), d_{nc}(-\lambda_j)).$$

(23)

Then $B_n \leq B_{1n} + B_{2n} + B_{3n}$ for terms $B_{in}$ defined in the following.

Using the fact that $n^{-1} 2\pi L_{n0}(2\lambda_j) \leq (2j)^{-1} \mathbb{1}_{\{j \leq \lfloor n/4 \rfloor\}} + (n-2j)^{-1} \mathbb{1}_{\{j > \lfloor n/4 \rfloor\}}$, by Lemma 1, and the fact that $\mathrm{cum}(d_{nc}(\lambda_j), d_{nc}(\lambda_j)) \leq C\lambda_j^{-\alpha}(\lambda_j^{-1} + \log(n) L_{n0}(2\lambda_j))$, by Lemma 3(i), we have, for $B_{1n} = n^{-3} \sum_{j=1}^N |g(\lambda_j) h(\lambda_j)| \times \mathrm{cum}^2(d_{nc}(\lambda_j), d_{nc}(\lambda_j))$, that

$$B_{1n} \leq Cn^{-1+\max\{0, 2\alpha-2\beta\}} \log^2(n) \left( \sum_{j=1}^{\lfloor n/4 \rfloor} j^{-2} + \sum_{j=\lfloor n/4 \rfloor+1}^N (n-2j)^{-2} \right) = o(1).$$

By Lemma 3(ii), $|\mathrm{cum}(d_{nc}(\lambda_j), d_{nc}(\lambda_j), d_{nc}(-\lambda_j), d_{nc}(-\lambda_j))| \leq Cn(n^{1/2} + \log^3(n))\lambda_j^{-2\alpha}$ so that $B_{2n} = n^{-3} \sum_{j=1}^N |g(\lambda_j) h(\lambda_j) \mathrm{cum}(d_{nc}(\lambda_j), d_{nc}(\lambda_j),$

$d_{nc}(-\lambda_j), d_{nc}(-\lambda_j))| = o(1)$. Choose $0 < \rho < \pi$. Using Lemma 3(i) and Lemma 2, for an arbitrarily small $\rho$, we may bound

$$\overline{\lim} B_{3n} \leq \overline{\lim}\left(\frac{C}{n}\sum_{\lambda_1 \leq \lambda_j < \rho} \lambda_j^{2\beta-2\alpha} + \frac{R_{n\rho}}{n} \cdot \frac{C}{n}\sum_{\rho \leq \lambda_j \leq \lambda_N} \lambda_j^{2\beta-\alpha}\right)$$

$$\leq C \int_0^\rho \lambda^{2\beta-2\alpha}\, d\lambda$$

for $B_{3n} = |4\pi/n \sum_{j=1}^N g(\lambda_j) h(\lambda_j)[\text{cum}_j - f(\lambda_j)][\text{cum}_j + f(\lambda_j)]|$, where we denote $\text{cum}_j = \text{cum}(d_{nc}(\lambda_j), d_{nc}(-\lambda_j))/(2\pi n)$. As $C$ does not depend on $\rho$ above, $B_{3n} = o(1)$ follows. $\square$

LEMMA 8. *Suppose Assumption* A.1 *holds and* $0 \leq \beta < 1$, $\alpha - \beta < 1/2$. *Let* $\omega = \max\{1/3, 1/4 + (\alpha - \beta)/2\}$. *Then* $\max_{1 \leq j \leq N} I_n(\lambda_j)\lambda_j^\beta = o_p(n^\omega)$ *and* $\max_{1 \leq j \leq N} f(\lambda_j)\lambda_j^\beta = o(n^\omega)$.

PROOF. It holds that $\mathrm{E}(I_n^4(\lambda_j)) = \text{cum}(I_{nc}^2(\lambda_j), I_{nc}^2(\lambda_j)) + [\mathrm{E}(I_{nc}^2(\lambda_j))]^2 \leq C\lambda_j^{-4\alpha}$, by the product theorem for cumulants; see (33), (23) and Lemma 3. For each $\varepsilon > 0$, we then have

$$P\Big(\max_{1 \leq j \leq N} I_n(\lambda_j)\lambda_j^\beta > \varepsilon n^\omega\Big) \leq \frac{1}{\varepsilon n^\omega}\left(\sum_{j=1}^N \lambda_j^{4\beta}\mathrm{E}(I_{nc}^4(\lambda_j))\right)^{1/4} = o(1),$$

which follows from $n^{-4\omega}\sum_{j=1}^N \lambda_j^{4\beta-4\alpha} = o(1)$. $\square$

PROOF OF THEOREM 1. We give a detailed argument for Theorem 1(i); parts (ii)–(iii) of Theorem 1 follow with some minor modifications. By Lemma 5, $0 \in \text{ch}^\circ\{\pi G_{\theta_0}(\lambda_j)I_n(\lambda_j)\}_{j=1}^N \subset \mathbb{R}^r$ with probability approaching 1 as $n \to \infty$ so that a positive $R_n(\theta_0)$ exists in probability. In view of (15), we can express the extrema $R_n(\theta_0) = \prod_{j=1}^N (1+\gamma_j)^{-1}$ with $\gamma_j = t'_{\theta_0}\pi G_{\theta_0}(\lambda_j)I_n(\lambda_j)$, $|\gamma_j| < 1$, where $t_{\theta_0} \in \mathbb{R}^r$ satisfies $Q_{1n}(\theta_0, t_{\theta_0}) = 0$ for the function $Q_{1n}(\cdot,\cdot)$ on $\Theta \times \mathbb{R}^r$ defined in (28). Let

$$W_{n\theta_0} = \frac{2\pi}{n}\sum_{j=1}^N G_{\theta_0}(\lambda_j)G'_{\theta_0}(\lambda_j)I_n^2(\lambda_j),$$

(24)

$$J_{n\theta_0} = Q_{1n}(\theta_0, 0) = \frac{2\pi}{n}\sum_{j=1}^N G_{\theta_0}(\lambda_j)I_n(\lambda_j).$$

By Lemma 6 with $\int_\Pi G_{\theta_0} f\, d\lambda = \mathcal{M} = 0$ and Lemma 7, we have

(25) $\qquad |J_{n\theta_0}| = O_p(n^{-1/2}), \qquad \|W_{n\theta_0} - W_{\theta_0}\| = o_p(1)$

so that $W_{n\theta_0}$ is nonsingular in probability. Using Assumption A.2 and Lemma 8, it holds that

$$(26) \quad Y_n = \max_{1 \leq j \leq N} \pi \|G_{\theta_0}(\lambda_j)\| I_n(\lambda_j) = o_p(n^{1/2}), \qquad \|t_{\theta_0}\| = O_p(n^{-1/2}),$$

where the order of $\|t_{\theta_0}\|$ follows as in [33, 34]. Note that by (26), $\max_{1 \leq j \leq N} |\gamma_j| \leq \|t_{\theta_0}\| Y_n = O_p(n^{-1/2}) o_p(n^{1/2}) = o_p(1)$ holds. Algebraically, we write $0 = Q_{1n}(\theta_0, t_{\theta_0}) = J_{n\theta_0} - \pi W_{n\theta_0} t_{\theta_0} + (2\pi/n) \sum_{j=1}^{N} G_{\theta_0}(\lambda_j) I_n(\lambda_j) \gamma_j^2 / (1 + \gamma_j)$ and solve for $t_{\theta_0} = (\pi W_{n\theta_0})^{-1} J_{n\theta_0} + \phi_n$, where

$$\|\phi_n\| \leq Y_n \|t_{\theta_0}\|^2 \|W_{n\theta_0}^{-1}\| \left( \frac{2\pi}{n} \sum_{j=1}^{N} \|G_{\theta_0}(\lambda_j)\|^2 I_n^2(\lambda_j) \right) \left\{ \max_{1 \leq j \leq N} (1 + \gamma_j)^{-1} \right\}$$

$$= o_p(n^{-1/2}),$$

by Lemma 7, (26) and the fact that $\max_{1 \leq j \leq N} |\gamma_j| = o_p(1)$. When $\|t_{\theta_0}\| Y_n < 1$, we apply a Taylor expansion $\log(1 + \gamma_j) = \gamma_j - \gamma_j^2/2 + \Delta_j$ for each $1 \leq j \leq N$. Then

$$\ell_n(\theta_0) = 4 \sum_{j=1}^{N} \log(1 + \gamma_j) = 2 \left[ 2 \sum_{j=1}^{N} \gamma_j - \sum_{j=1}^{N} \gamma_j^2 \right] + 4 \sum_{j=1}^{N} \Delta_j,$$

(27)

$$2 \left[ 2 \sum_{j=1}^{N} \gamma_i - \sum_{j=1}^{N} \gamma_j^2 \right] = n J'_{n\theta_0} (\pi W_{n\theta_0})^{-1} J_{n\theta_0} - n \phi'_n (\pi W_{n\theta_0}) \phi_n.$$

By Lemma 6 and (25), $n J'_{n\theta_0} (\pi W_{n\theta_0})^{-1} J_{n\theta_0} \xrightarrow{d} \chi_r^2$. We also have $n \phi'_n (\pi W_{n\theta_0}) \times \phi_n = o_p(1)$ and we may bound $\sum_{j=1}^{N} |\Delta_j|$ by

$$\frac{n \|t_{\theta_0}\|^3 Y_n}{(1 - \|t_{\theta_0}\| Y_n)^3} \frac{1}{n} \sum_{j=1}^{N} \|\pi G_{\theta_0}(\lambda_j)\|^2 I_n^2(\lambda_j) = n O_p(n^{-3/2}) o_p(n^{1/2}) O_p(1)$$

$$= o_p(1),$$

from Lemma 7 and (26). Theorem 1(i) follows by Slutsky's Theorem.

Under Theorem 1(ii)–(iii), $f_{\theta_0} = f$ holds and $R_{n,\mathcal{F}}(\theta_0)$ exists in probability as $n \to \infty$, by Assumption A.5. We repeat the same arguments as above, replacing each occurrence of $I_n(\lambda_j)$ with $I_n(\lambda_j) - f(\lambda_j)$ instead; we denote the resulting quantities with a tilde:

$$\tilde{W}_{n\theta_0} = \frac{2\pi}{n} \sum_{j=1}^{N} G_{\theta_0}(\lambda_j) G'_{\theta_0}(\lambda_j) (I_n(\lambda_j) - f(\lambda_j))^2,$$

$$\tilde{J}_{n\theta_0} = \frac{2\pi}{n} \sum_{j=1}^{N} G_{\theta_0}(\lambda_j) (I_n(\lambda_j) - f(\lambda_j)),$$

$\tilde{\gamma}_j$, $\tilde{\Delta}_j$, $\tilde{\phi}_n$, etc. All the previous points follow except for two, which are straightforward to remedy: by Lemma 7, $\|2\tilde{W}_{n\theta_0} - W_{\theta_0}\| = o_p(1)$ instead of (25), and in (27), we must write

$$\ell_{n,\mathcal{F}}(\theta_0) = 2\sum_{j=1}^{N} \log(1+\tilde{\gamma}_j)$$

$$= n\tilde{J}'_{n\theta_0}(2\pi\tilde{W}_{n\theta_0})^{-1}\tilde{J}_{n\theta_0} - n\tilde{\phi}'_n(2^{-1}\pi\tilde{W}_{n\theta_0})\tilde{\phi}_n + 2\sum_{j=1}^{N}\tilde{\Delta}_j,$$

where $n\tilde{J}'_{n\theta_0}(2\pi\tilde{W}_{n\theta_0})^{-1}\tilde{J}_{n\theta_0} \xrightarrow{d} \chi_r^2$ by Lemma 6, $n\tilde{\phi}'_n(2^{-1}\pi\tilde{W}_{n\theta_0})\tilde{\phi}_n = o_p(1)$ and $\sum_{j=1}^{N}\tilde{\Delta}_j = o_p(1)$. □

PROOF OF THEOREM 2. We require some additional notation. Define the following functions on $\Theta \times \mathbb{R}^r$:

(28)
$$Q_{1n}(\theta,t) = \frac{2\pi}{n}\sum_{j=1}^{N}\frac{G_\theta(\lambda_j)I_n(\lambda_j)}{1+t'\pi G_\theta(\lambda_j)I_n(\lambda_j)},$$

$$Q_{2n}(\theta,t) = \frac{2\pi}{n}\sum_{j=1}^{N}\frac{I_n(\lambda_j)(\partial G_\theta(\lambda_j)/\partial\theta)'t}{1+t'\pi G_\theta(\lambda_j)I_n(\lambda_j)}.$$

Also, define versions $\tilde{Q}_{1n}(\theta,t)$ and $\tilde{Q}_{2n}(\theta,t)$ by replacing each $I_n(\lambda_i)$ with $I_n(\lambda_i) - f_\theta(\lambda_i)$ in (28) and adding the extra term $-(\partial f_\theta(\lambda_i)/\partial\theta)G'_\theta(\lambda_i)t$ to the numerator of $Q_{2n}$. We use the following MELE existence result to prove Theorem 2; see [32] for its proof:

LEMMA 9. *Under the assumptions of Theorem 2,*

*(i) the probability that $R_n(\theta)$ attains a maximum $\hat{\theta}_n$, which lies in the ball $\mathcal{B}_n$ and satisfies $Q_{1n}(\hat{\theta}_n, t_{\hat{\theta}_n}) = 0$ and $\partial\ell_n(\theta)/\partial\theta|_{\theta=\hat{\theta}_n} = 2nQ_{2n}(\hat{\theta}_n, t_{\hat{\theta}_n}) = 0$, converges to 1 as $n \to \infty$;*

*(ii) under the assumptions of Theorem 2(ii) or (iii), result (i) above holds for $R_{n,\mathcal{F}}(\theta)$ with respect to $\hat{\theta}_{n,\mathcal{F}}$, $\ell_{n,\mathcal{F}}$, $\tilde{Q}_{1n}(\hat{\theta}_{n,\mathcal{F}}, t_{\hat{\theta}_{n,\mathcal{F}}})$, $2^{-1}\tilde{Q}_{2n}(\hat{\theta}_{n,\mathcal{F}}, t_{\hat{\theta}_{n,\mathcal{F}}})$ [replacing $\hat{\theta}_n$, $\ell_n$, $Q_{1n}(\hat{\theta}_n, t_{\hat{\theta}_n})$, $Q_{2n}(\hat{\theta}_n, t_{\hat{\theta}_n})$].*

We now establish the asymptotic normality of $\hat{\theta}_n$, following arguments in [38]. Expanding $Q_{1n}(\hat{\theta}_n, t_{\hat{\theta}_n})$, $Q_{2n}(\hat{\theta}_n, t_{\hat{\theta}_n})$ at $(\theta_0, 0)$ with Lemma 9, we have that

$$\Sigma_n\begin{pmatrix} t_{\hat{\theta}_n} \\ \hat{\theta}_n - \theta_0 \end{pmatrix} = \begin{bmatrix} -J_{n\theta_0} + E_{1n} \\ E_{2n} \end{bmatrix},$$

$$\Sigma_n = \begin{bmatrix} \partial Q_{1n}(\theta_0,0)/\partial t & \partial Q_{1n}(\theta_0,0)/\partial \theta \\ \partial Q_{2n}(\theta_0,0)/\partial t & 0 \end{bmatrix},$$

$$Q_{1n}(\theta_0,0) = J_{n\theta_0}, \qquad \frac{\partial Q_{1n}(\theta_0,0)}{\partial t} = -\pi W_{n\theta_0},$$

$$\frac{\partial Q_{2n}(\theta_0,0)}{\partial t} = \left[\frac{\partial Q_{1n}(\theta_0,0)}{\partial \theta}\right]', \qquad \frac{\partial Q_{2n}(\theta_0,0)}{\partial \theta} = 0,$$

with $J_{n\theta_0}, W_{n\theta_0}$ as in (24). In addition, it can be shown that

$$\frac{\partial Q_{1n}(\theta_0,0)}{\partial \theta} = \frac{2\pi}{n} \sum_{j=1}^{N} \frac{\partial G_{\theta_0}(\lambda_j)}{\partial \theta} I_n(\lambda_j) \xrightarrow{p} \frac{1}{2} \int_{\Pi} \frac{\partial G_{\theta_0}(\lambda)}{\partial \theta} f(\lambda)\, d\lambda = \frac{D_{\theta_0}}{2},$$

using $(2\pi/n) \sum_{j=1}^{N} f(\lambda_j) \partial G_{\theta_0}(\lambda_j)/\partial \theta \to D_{\theta_0}/2$ by the dominated convergence theorem and a modification of the proof of Lemma 7. Applying this convergence result with (25),

$$(29) \qquad \Sigma_n^{-1} = \begin{bmatrix} A_{11n} & A_{12n} \\ A_{21n} & A_{22n} \end{bmatrix} \xrightarrow{p} \frac{1}{2\pi} \begin{bmatrix} -2\pi U_{\theta_0} & W_{\theta_0}^{-1} D_{\theta_0} V_{\theta_0} \\ V_{\theta_0} D'_{\theta_0} W_{\theta_0}^{-1} & 2\pi V_{\theta_0} \end{bmatrix}$$

holds. One may verify that $\|E_{1n}\|$, $\|E_{2n}\| = O_p(\delta_n n^{-1-\eta} \sum_{j=1}^{N} \Psi_{jn} \times \max\{I_n^2(\lambda_j), I_n(\lambda_j)\}) = o_p(\delta_n)$ for $\delta_n = \|t_{\hat{\theta}_n}\| + \|\hat{\theta}_n - \theta_0\|$ and $\Psi_{jn} = n^{-\eta} \lambda_j^{-2\delta} + \lambda_j^{\beta-\delta} + \lambda_j^{2\beta}$. By Lemma 6, $\sqrt{n} J_{n\theta_0} \xrightarrow{d} \mathcal{N}(0, \pi W_{\theta_0})$, and so it follows that $\delta_n = O_p(n^{-1/2})$. We then have that

$$(30) \qquad \begin{aligned} \sqrt{n}(\hat{\theta}_n - \theta_0) &= -\sqrt{n} A_{21n} J_{n\theta_0} + o_p(1) \xrightarrow{d} \mathcal{N}(0, V_{\theta_0}), \\ \sqrt{n}(t_{\hat{\theta}_n} - 0) &= -\sqrt{n} A_{11n} J_{n\theta_0} + o_p(1) \xrightarrow{d} \mathcal{N}(0, U_{\theta_0}). \end{aligned}$$

To show the normality of $\hat{\theta}_{n,\mathcal{F}} - \theta_0$ under the conditions of Theorem 2(ii)–(ii), we substitute $\tilde{Q}_{1n}, \tilde{Q}_{2n}$ for $Q_{1n}, Q_{2n}$ above and, using Lemma 9, repeat the same expansion with an analogously defined matrix $\tilde{\Sigma}_n$ (having components $\tilde{A}_{ijn}$ in $\tilde{\Sigma}_n^{-1}$). Note that $\partial \tilde{Q}_{1n}(\theta_0,0)/\partial \theta = \partial Q_{1n}(\theta_0,0)/\partial \theta - \tilde{D}_{n\theta_0}$ where, by the Lebesgue dominated convergence theorem, as $n \to \infty$,

$$\tilde{D}_{n\theta_0} = \frac{2\pi}{n} \sum_{i=1}^{N} \left.\frac{\partial [G_\theta(\lambda_i) f_\theta(\lambda_i)]}{\partial \theta}\right|_{\theta=\theta_0} \to \frac{1}{2} \frac{\partial}{\partial \theta} \left[\int_{\Pi} f_\theta G_\theta\, d\lambda\right]\bigg|_{\theta=\theta_0} = 0,$$

since the theorem's conditions justify exchanging the order of differentiation/integration of $f_\theta G_\theta$ at $\theta_0$ and $\int_\Pi f_\theta G_\theta\, d\lambda = \mathcal{M}$ is constant in a neighborhood of $\theta_0$. The convergence result in (29) follows for $\tilde{\Sigma}_n^{-1}$ upon replacing "$-2\pi U_{\theta_0}, 2\pi V_{\theta_0}$" with "$-4\pi U_{\theta_0}, \pi V_{\theta_0}$." By Lemma 6, $\sqrt{n}\tilde{J}_{n\theta_0} = \sqrt{n}\tilde{Q}_{1n}(\theta_0,0) \xrightarrow{d} \mathcal{N}(0, \pi W_{\theta_0})$ so that (30) holds for $\sqrt{n}(\hat{\theta}_{n,\mathcal{F}} - \theta_0)$ and $\sqrt{n}(t_{\hat{\theta}_{n,\mathcal{F}}} - 0)/2$ after replacing $J_{n\theta_0}, A_{21n}, A_{11n}$ with $\tilde{J}_{n\theta_0}, \tilde{A}_{21n}, \tilde{A}_{11n}/2$. $\square$

PROOF OF THEOREM 3. To establish (i), assume without loss of generality that $G_\theta(\lambda)$ is real-valued and increasing in $\theta$ (i.e., $r = 1$). For $\theta \in \Theta$, define $E_{n,\theta} = \frac{1}{N} \sum_{j=1}^{N} \{\pi I_n(\lambda_j) G_\theta(\lambda_j) - \mathcal{M}\}$, where $\mathcal{M} = \int_0^\pi G_{\theta_0} f \, d\lambda$. Following the proof of Lemma 7, $E_{n,\theta} \xrightarrow{p} \int_0^\pi G_\theta f \, d\lambda - \mathcal{M}$ holds for each $\theta$ so that, for an arbitrarily small $\varepsilon > 0$, $P(E_{n,\theta_0-\varepsilon} < 0 < E_{n,\theta_0+\varepsilon}) \to 1$ by the monotonicity of $G_\theta$. When the event in the probability statement holds, there exists $\hat{\theta}_n \in \Theta$ where $E_{n,\hat{\theta}_n} = 0$, by the continuity of $G_\theta$, and $L_n(\hat{\theta}_n) = (\pi/N)^N$ follows in (13). For any $\theta$ with $|\theta - \theta_0| \geq \varepsilon$, we have that $E_{n,\theta} \neq 0$, implying that $L_n(\theta) < L_n(\hat{\theta}_n)$. Hence, a global maximum $\hat{\theta}_n$ satisfies $P(|\hat{\theta}_n - \theta_0| < \varepsilon) \to 1$.

For (ii), we consider $G_\theta \equiv G_\theta^w$ and $\mathcal{M} \equiv \mathcal{M}_w$ from (11). Suppose $\theta_n^*$ and $\theta^*$ denote the minimums of $W_n(\theta) = \frac{1}{4} \log \frac{\sigma^2}{2\pi} + \frac{1}{4\pi N} \sum_{j=1}^N \pi I_n(\lambda_j) f_\theta^{-1}(\lambda_j)$ and $W(\theta)$, respectively. Using Lemma 4 with arguments as in Theorem 1 of [20], it follows that $\theta_n^* \xrightarrow{p} \theta^*$. Since $\theta^*$ is interior to $\Theta$, it holds that $\partial W_n(\theta_n^*)/\partial \theta = 0$, which implies that $E_{n,\theta_n^*} = 0$ under the above definition of $E_{n,\theta}$ with $G_\theta^w, \mathcal{M}_w$. Hence, $L_n(\theta_n^*) = (\pi/N)^N$ and, using Lemma 4 with arguments of Lemma 1 of [20], $0 = E_{n,\theta_n^*} \xrightarrow{p} \int_0^\pi G_{\theta^*} f \, d\lambda - \mathcal{M}$ holds, whereby $\theta^* = \theta_0$, by uniqueness. We then have a global maximum $\hat{\theta}_n = \theta_n^*$, for which $\hat{\theta}_n \xrightarrow{p} \theta_0$. □

PROOF OF THEOREM 4. Considering Theorem 4(i), let $P_X = X(X'X)^{-1}X'$ and $I_{r \times r}$ denote the projection matrix for a given matrix $X$ and the $r \times r$ identity matrix. Writing $(\pi W_{\theta_0}/n)^{1/2} Z_{n\theta} = J_{n\theta_0} + D_{\theta_0}(\theta - \theta_0)/2$, it holds that $|\ell(\theta) - Z'_{n\theta} Z_{n\theta}| = o_p(1 + n\|\theta - \theta_0\|^2)$ uniformly for $\theta \in \mathcal{B}_n$ (see [32], Theorem 3) so that $\ell_n(\hat{\theta}_n) = Z'_{n\theta_0}(I_{r \times r} - P_{W_{\theta_0}^{-1/2} D_{\theta_0}}) Z_{n\theta_0} + o_p(1)$, using (29)–(30). By (27), $\ell_n(\theta_0) = Z'_{n\theta_0} Z_{n\theta_0} + o_p(1)$, where $Z_{n\theta_0} \xrightarrow{d} \mathcal{N}(0, I_{r \times r})$, by Lemma 6. Theorem 4 follows since $P_{W_{\theta_0}^{-1/2} D_{\theta_0}}$, $I_{r \times r} - P_{W_{\theta_0}^{-1/2} D_{\theta_0}}$ are orthogonal idempotent matrices with ranks $r$, $r - p$. □

PROOFS OF COROLLARIES 1 AND 2. Nordman [32] gives a detailed proof of Corollary 1 and Corollary 2 follows by modifying arguments from Corollary 5 of [38]. □

PROOF OF THEOREM 5. Nordman [32] provides details where the most important, additional distributional results required are $(2\pi/\sqrt{n}) \sum_{j=1}^N Y_{n\theta_0,j} \xrightarrow{d} \mathcal{N}(0, \pi W_{\theta_0})$, $(2\pi/n) \sum_{j=1}^N Y_{n\theta_0,j} Y'_{n\theta_0,j} \xrightarrow{p} W_{\theta_0}/2$ for $Y_{n\theta_0,j} = (I_{nc}^2(\lambda_j)/[2f_{\theta_0}^2(\lambda_j)] - \pi, I_{nc}(\lambda_j) G_{\theta_0}^w(\lambda_j) - \mathcal{M}_w)'$ and

$$W_{\theta_0} = \begin{bmatrix} W_{\theta_0}^* & 0 \\ 0 & W_{\theta_0}^{**} \end{bmatrix}, \qquad W_{\theta_0}^* = \begin{bmatrix} 10\pi & 4\pi \\ 4\pi & 2\pi \end{bmatrix},$$

$$W^{**}_{\theta_0} = \left(\int_\Pi f^2_{\theta_0} \frac{\partial f^{-1}_{\theta_0}}{\partial \vartheta_i} \frac{\partial f^{-1}_{\theta_0}}{\partial \vartheta_j} d\lambda\right)_{i,j=1,\ldots,p-1}.$$

The convergence results can be shown using arguments in [3] and [16]. □

## APPENDIX

LEMMA A.1. *Suppose Assumptions* A.1–A.2 *hold for a real-valued, even function $g$ satisfying Condition 2 of Assumption A.3. Define the nth Cesàro mean of the Fourier series of $g$ as $c_n g(\lambda) = \int_\Pi K_n(\lambda - y) g(y)\, dy$, $\lambda \in \Pi$, and let $\mathcal{C}_n = \int_0^\pi c_n g(\lambda) I_{nc}(\lambda)\, d\lambda$. Then as $n \to \infty$,*

(31)
$$\sqrt{n}\left|\mathcal{C}_n - \int_0^\pi g(\lambda) I_{nc}(\lambda)\, d\lambda\right| = o_p(1),$$

$$\sqrt{n}\left|\frac{2\pi}{n}\sum_{j=1}^N g(\lambda_j) I_n(\lambda_j) - \mathcal{C}_n\right| = o_p(1).$$

PROOF. Assume that $g$ satisfies A.1 and Condition 2 with respect to $\beta$ and $\beta_1$, respectively. We establish the first part of (31). Using $\int_\Pi K_n(\lambda - y)\, dy = 1$, $\mathrm{E}(I_{nc}(\lambda)) = \int_\Pi K_n(\lambda - y) f(y)\, dy$ and Fubini's Theorem, we write $2\mathrm{E}|\int_0^\pi (c_n g - g) I_{nc}\, d\lambda|$ as

$$\mathrm{E}\left|\int_\Pi I_{nc}(\lambda)\left[\int_\Pi K_n(\lambda - y)[g(y) - g(\lambda)]\, dy\right] d\lambda\right|$$

$$\leq \int_{\Pi^3} K_n(\lambda - z) K_n(\lambda - y) f(z) |g(y) - g(\lambda)|\, dz\, dy\, d\lambda \leq t_{1n} + t_{2n},$$

$t_{1n} = \int_{\Pi^2} K_n(\lambda - z) f(z) |g(z) - g(\lambda)|\, dz\, d\lambda$, $t_{2n} = \int_{\Pi^3} K_n(\lambda - z) K_n(y - \lambda) f(z) \times |g(y) - g(z)|\, d\lambda\, dz\, dy$. It suffices to show that $\sqrt{n} t_{2n} = o(1)$ since arguments from (22) provide $\sqrt{n} t_{1n} = o(1)$. Applying Lemma 1(iv) and (20) sequentially, we bound $\int_\Pi K_n(\lambda - z) K_n(y - \lambda)\, d\lambda$ by

$$\frac{C}{n^2}\int_\Pi L^2_{n0}(\lambda - z) L^2_{n0}(y - \lambda)\, d\lambda \leq \frac{Cn}{(1 + |(y-z), \mathrm{mod}\, 2\pi| n)^2}.$$

From this and arguments from (22), we have

$$\sqrt{n} t_{2n} \leq C n^{3/2} \int_{\Pi^2} \frac{f(z)|g(y) - g(z)|}{(1 + |(y-z)\, \mathrm{mod}\, 2\pi| n)^2}\, dz\, dy = o(1).$$

For the second part of (31), we use $(2\pi/n)\sum_{j=-N}^{\lfloor n/2 \rfloor} c_n g(\lambda_j) I_{nc}(\lambda_j) = 2\mathcal{C}_n$ and write

$$2\sqrt{n}\left|\frac{2\pi}{n}\sum_{j=1}^N g(\lambda_j) I_n(\lambda_j) - \mathcal{C}_n\right| \leq 4\pi\sqrt{n}(t_{3n} + t_{4n}),$$

$t_{3n} = n^{-1}\sum_{j=1}^{N} I_n(\lambda_j)|c_n g(\lambda_j) - g(\lambda_j)|$, $t_{4n} = n^{-1}(|c_n g(0)|I_{nc}(0) + |c_n g(\pi)| \times I_{nc}(\pi))$. We have $\sqrt{n}t_{4n} = o_p(1)$ from $\mathrm{E}[I_{nc}(0)] \leq Cn^\alpha$, $|c_n g(0)| \leq Cn^{-\beta}$ and $|c_n g(\pi)|\mathrm{E}[I_{nc}(\pi)] \leq C$ by Assumptions A.1–A.2. For a $C > 0$ independent of $1 \leq j \leq N$ ($n > 3$), if we establish

$$(32) \quad |c_n g(\lambda_j) - g(\lambda_j)| \leq C\lambda_j^{\beta_1}\left(\frac{\log(n)}{j} + \frac{\mathbb{1}_{\{j>n/4\}}}{n-2j}\right), \qquad 1 \leq j \leq N,$$

then it will follow that $\sqrt{n}t_{3n} \leq Cn^{-1/2}\log(n)(\max_{1\leq j\leq N} I_n(\lambda_j)\lambda_j^{\beta_1})\sum_{j=1}^{n} j^{-1} = o_p(1)$ from Lemma 8.

Fix $1 \leq j \leq N$. To prove (32), we decompose the difference

$$|c_n g(\lambda_j) - g(\lambda_j)| = \left|\int_\Pi K_n(y)[g(\lambda_j - y) - g(\lambda_j)]\,dy\right|$$

$$\leq \left|\int_0^\pi d_{jn}^+(y)\,dy\right| + \left|\int_0^\pi d_{jn}^-(y)\,dy\right|,$$

where $d_{jn}^+(y) = K_n(y)[g(\lambda_j + y) - g(\lambda_j)]$, $d_{jn}^-(y) = K_n(y)[g(\lambda_j - y) - g(\lambda_j)]$. We separately bound the last two absolute integrals using (20), $|g(y) - g(z)| \leq C\,|y|^{-1+\beta_1}|z - y|$ for $0 < |y| \leq |z| \leq \pi$ and the fact that $\lambda_j > \pi/2$ if and only if $j > n/4$. Using Lemma 1, we bound

$$\left|\int_0^{1/n} d_{jn}^+(y)\,dy\right| \leq C\lambda_j^{-1+\beta_1}\int_0^{1/n} L_{n0}(y)y\,dy \leq Cj^{-1}\lambda_j^{\beta_1},$$

$$\left|\int_{1/n}^{\pi-\lambda_j} d_{jn}^+(y)\,dy\right| \leq C\lambda_j^{-1+\beta_1}\int_{1/n}^{\pi-\lambda_j} \frac{ny}{(1+ny)^2}\,dy \leq Cj^{-1}\log(n)\lambda_j^{\beta_1},$$

and if $\lambda_j \leq \pi/2$,

$$\left|\int_{\pi-\lambda_j}^\pi d_{jn}^+(y)\,dy\right| \leq Cn^{-1}L_0^2(\pi/2)\int_{\pi-\lambda_j}^\pi 1\,dy \leq Cn^{-1}\lambda_j.$$

If $\lambda_j > \pi/2$, we use a substitution $u = 2\pi - \lambda_j - y$ and the fact that $1 \leq n - 2j \leq N$ to find

$$\left|\int_{\pi-\lambda_j}^{2\pi-2\lambda_j} d_{jn}^+(y)\,dy\right| = \left|\int_{\lambda_j}^\pi K_n(2\pi - u - \lambda_j)[g(u) - g(\lambda_j)]\,du\right|$$

$$\leq C\lambda_j^{\beta_1}\int_{\lambda_j}^\pi \frac{(u-\lambda_j)n}{(1+(2\pi - u - \lambda_j)n)^2}\,du \leq Cn^{-1}\log(n)\lambda_j^{\beta_1},$$

$$\left|\int_{2\pi-2\lambda_j}^\pi d_{jn}^+(y)\,dy\right| \leq C\int_{\pi-\lambda_j}^{\lambda_j} K_n(2\pi - u - \lambda_j)\,du \leq C(n-2j)^{-1}\lambda_j^{\beta_1}.$$

Hence, the bound in (32) applies to $|\int_0^\pi d_{jn}^+\,dy|$; the same can be shown for $|\int_0^\pi d_{jn}^-\,dy|$ by considering separate integrals over the intervals $(0, 1/n]$, $(1/n, \lambda_j/2]$, $(\lambda_j/2, \lambda_j]$ and $(\lambda_j, \pi]$ if $\lambda_j > \pi/2$, or $(\lambda_j, 2\lambda_j]$ and $(2\lambda_j, \pi]$ if $\lambda_j \leq \pi/2$. See Lemma 12 of [32]. □

LEMMA A.2. *Under Assumption* A.1, *suppose* $g, h$ *are real-valued, even Riemann integrable functions on* $\Pi$ *such that* $|g(\lambda)|, |h(\lambda)| \leq C|\lambda|^\beta$, $0 \leq \beta < 1$, $\alpha - \beta < 1/2$. *Then as* $n \to \infty$, $V_n = \text{Var}[(2\pi/n) \sum_{j=1}^N g(\lambda_j) h(\lambda_j) I_n^2(\lambda_j)] = o(1)$.

PROOF. We expand $V_n$ and then bound $V_n \leq V_{1n} + V_{2n}$ with $V_{1n}$ and $V_{2n}$ defined as

$$\frac{C}{n^6} \sum_{j=1}^N (\lambda_j)^{4\beta} \text{cum}(|d_{nc}(\lambda_j)|^4, |d_{nc}(\lambda_j)|^4),$$

$$\frac{C}{n^6} \sum_{1 \leq i < j \leq N} (\lambda_i \lambda_j)^{2\beta} |\text{cum}(|d_{nc}(\lambda_i)|^4, |d_{nc}(\lambda_j)|^4)|,$$

respectively. Let $\mathcal{P}$ be the set of all indecomposable partitions $P$ of the labels in the two row table $\{a_{st}\}$, $s = 1, 2$, $t = 1, 2, 3, 4$ ([5], Section 2.3). We write the elements of a partition $P = (P_1, \ldots, P_m)$, $1 \leq m \leq 7$, with each $P_i \subset \{a_{st}\} = \bigcup_{i=1}^m P_i$, $P_i \cap P_j = \varnothing$ if $i \neq j$. For $1 \leq i \leq j \leq N$, we define $a_{11}^{ij} = a_{12}^{ij} = -a_{13}^{ij} = -a_{14}^{ij} = \lambda_i$, $a_{21}^{ij} = a_{22}^{ij} = -a_{23}^{ij} = -a_{24}^{ij} = \lambda_j$. By the product theorem for cumulants,

$$\text{cum}(|d_n(\lambda_i)|^4, |d_n(\lambda_j)|^4) = \sum_{P \in \mathcal{P}} \text{cum}_{ijn}(P),$$

(33)

$$\text{cum}_{ijn}(P) = \prod_{u=1}^m \text{cum}(d_{nc}(a_{st}^{ij}) : a_{st} \in P_u).$$

Because $E(d_{nc}(\lambda)) = 0$, we need only consider those partitions $\mathcal{P}^* \equiv \{P = (P_1, \ldots, P_m) \in \mathcal{P} : 1 < |P_1| \leq \cdots \leq |P_m|, 1 \leq m \leq 6\}$ where each set in the partition has two or more elements $a_{st}$, using $|B|$ to denote the size of a finite set $B$. Defining $U_{1n}(P) = n^{-6} \sum_{j=1}^N (\lambda_j)^{4\beta} |\text{cum}_{jjn}(P)|$ and $U_{2n}(P) = n^{-6} \sum_{1 \leq i < j \leq N} (\lambda_i \lambda_j)^{2\beta} |\text{cum}_{ijn}(P)|$, we can bound $V_{in} \leq C \sum_{P \in \mathcal{P}^*} U_{in}(P)$, $i = 1, 2$, so that it suffices to show

(34) $\qquad U_{in}(P) = o(1), \qquad P \in \mathcal{P}^*, \ i = 1, 2.$

By Lemma 3, we have $|\text{cum}_{jjn}(P)| \leq Cn^4 \lambda_j^{-4\alpha}$ for $P \in \mathcal{P}^*, 1 \leq j \leq N$ so that $U_{1n}(P) \leq Cn^{\max\{0, 2\alpha - 2\beta\} - 1} (n^{-1} \sum_{j=1}^N \lambda_j^{2\beta - 2\alpha}) = o(1)$ since $\alpha - \beta < 1/2$. Hence, (34) is established for $U_{1n}$ and $V_{1n} = o(1)$.

We now show (34) for $U_{2n}(P)$ by bounding $|\text{cum}_{ijn}(P)|$, over $1 \leq i < j \leq N$, for several cases of $P = (P_1, \ldots, P_m) \in \mathcal{P}^*$. These cases are $m = 1$; $m = 2$, $|P_2| = 6$; $m = 2$, $|P_2| = 5$; $m = 2$, $|P_2| = 4$; $m = 3$, $|P_3| = |P_2| = 3$; $m = 3$, $|P_3| = 4$, $|P_2| = 2$. The last (seventh) case $m = 4$, $|P_1| = |P_2| = |P_3| = |P_4| = 2$ has the following subcases:

(7.1) there exist $k_1 \neq k_2$ where $\sum_{a_{st} \in P_{k_1}} a_{st}^{ij} = 0 = \sum_{a_{st} \in P_{k_2}} a_{st}^{ij}$;

(7.2) there exists exactly one $k$ where $\sum_{a_{st} \in P_k} a_{st}^{ij} = 0$;

(7.3) for each $k$, $\sum_{a_{st} \in P_k} a_{st}^{ij} \neq 0$ holds and for some $k_1, k_2$, $|\sum_{a_{st} \in P_{k_1}} a_{st}^{ij}| = 2\lambda_i$, $|\sum_{a_{st} \in P_{k_2}} a_{st}^{ij}| = 2\lambda_j$;

(7.4) for each $k$, $|\sum_{a_{st} \in P_k} a_{st}^{ij}| \notin \{0, 2\lambda_i, 2\lambda_j\}$.

The first six cases follow from Lemma 3 and (20). For example, under cases 3 or 4, we apply Lemma 3(ii) and (20) twice to bound $|\text{cum}_{ijn}(P)| \leq C\{n^{3/2} + n\log^4(n)\}^2 (\lambda_i \lambda_j)^{-2\alpha}$ and then $U_{2n}(P) \leq Cn^{-1}(n^{-1} \sum_{j=1}^{N} \lambda_j^{2\beta-2\alpha})^2 = o(1)$. See Lemma 13 of [32] for details.

To treat case 7, we define some sets. For $0 < \rho < \pi/2$, let $A = \{(i,j) : 1 \leq i < j \leq N\}$, $A_\rho = \{(i,j) \in A : \lambda_j < \rho\}$, $A^\rho = \{(i,j) \in A : \lambda_j \geq \rho\}$, $A^{n/2} = \{(i,j) \in A : i + j > n/2\}$ and $A_{n/2} = \{(i,j) \in A : i + j \leq n/2\}$. We will also use the fact that for integers $j > i \geq 1$,

$$(35) \quad \frac{1}{i(j-i)} \leq \frac{2}{j} \quad \text{if } i = 1 \text{ or } j(i-1) > i^2, \quad \text{otherwise,} \quad \frac{1}{i} < \frac{2}{j}.$$

For subcase 7.1, define $A_{\rho 1} = \{(i,j) \in A_\rho : i = 1 \text{ or } j(i-1) > i^2\}$ and $A_{\rho 2} = A_\rho \setminus A_{\rho 1}$.

*Subcase* 7.1. Without loss of generality, say $k_1 = 1$, $k_2 = 2$. We have $|\sum_{a_{st} \in P_3} a_{st}^{ij}| = |\sum_{a_{st} \in P_4} a_{st}^{ij}| \in \{\lambda_j - \lambda_i, \lambda_j + \lambda_i\}$ and by Lemma 3(i), $\prod_{k=1}^{2} \times |\text{cum}(d_{nc}(a_{st}^{ij}) : a_{st} \in P_k)| \leq Cn^2(\lambda_i \lambda_j)^{-\alpha}$. Fix $0 < \rho < \pi/2$. If $|\sum_{a_{st} \in P_k} a_{st}^{ij}| = \lambda_j - \lambda_i$, $k = 3, 4$, then by Lemma 1(i), Lemma 2 and (20) (for the sum involving $\lambda_j \geq \rho$) or Lemma 1(i), Lemma 3(i) and (20) (for the sum involving $\lambda_j < \rho$):

$$U_{2n}(P) \leq \frac{C(\rho)}{n^4} \sum_{A^\rho} \lambda_i^{\beta-\alpha} \left[\frac{n}{j-i} + R_{n\rho}\right]^2$$

$$+ \frac{C}{n^4} \sum_{A_\rho} \lambda_i^{2\beta-3\alpha} \lambda_j^{2\beta-\alpha} \left[\lambda_j^{-1} + \frac{n}{(j-i)^d}\right]^2$$

$$\equiv u_{1n}(\rho) + u_{2n}(\rho),$$

$$u_{1n}(\rho) \leq C(\rho)\left(n^{-1+\max\{0,\alpha-\beta\}} \sum_{j=1}^{n} j^{-2} + (n^{-1} R_{n\rho})^2 \left(\frac{1}{n} \sum_{j=1}^{N} \lambda_j^{\beta-\alpha}\right)\right) = o(1),$$

with some fixed $\max\{\alpha, 1/2\} < d < 1$. Using (35) on the sums over $A_{\rho 1}$ and $A_{\rho 2}$, we have

$$u_{2n}(\rho) \leq C\left(\frac{1}{n} \sum_{\lambda_1 \leq \lambda_i < \rho} \lambda_i^{2\beta-2\alpha}\right)^2 + \frac{C}{n^2} \sum_{A_{\rho 1}} \frac{(\lambda_i \lambda_j)^{2\beta-2\alpha}}{(j-i)^{d-\alpha}} + \frac{C}{n^2} \sum_{A_{\rho 2}} \frac{(\lambda_i \lambda_j)^{2\beta-2\alpha}}{(j-i)^d}$$

$$\leq C\bigg(\frac{1}{n}\sum_{\lambda_1\leq\lambda_i<\rho}\lambda_i^{2\beta-2\alpha}\bigg)^2.$$

Then $\overline{\lim}\,(u_{1n}(\rho)+u_{2n}(\rho))\leq C(\int_0^\rho \lambda^{2\beta-2\alpha}\,d\lambda)^2 = C\rho^{2+4\beta-4\alpha}$ for a $C$ that does not depend on $0<\rho<\pi/2$. Hence, $U_{2n}(P)=o(1)$ since $2+4\beta-4\alpha>0$. If $|\sum_{a_{st}\in P_k} a_{st}^{ij}|=\lambda_i+\lambda_j$, $k=3,4$, we essentially repeat the same steps as above to show $U_{2n}(P)\leq u_{3n}(\rho)+u_{4n}(\rho)$, where $u_{4n}(\rho)\leq Cu_{2n}(\rho)$, and, using Lemma 1(i) to expand $u_{3n}(\rho)=C(\rho)/n^4 \sum_{A^\rho} \lambda_i^{\beta-\alpha}[L_{n0}(\lambda_i+\lambda_j)+R_{n\rho}]^2$,

$$u_{3n}(\rho)\leq \frac{C(\rho)}{n^2}\bigg(\sum_{A^{n/2}}\frac{\lambda_i^{\beta-\alpha}}{(n-j-i)^2}+\sum_{A_{n/2}}\frac{\lambda_i^{\beta-\alpha}}{(j+i)^2}+(R_{n\rho})^2\bigg(\frac{1}{n}\sum_{j=1}^N \lambda_i^{\beta-\alpha}\bigg)\bigg)$$

$$\leq C(\rho)\bigg(n^{-1+\max\{0,\alpha-\beta\}}\sum_{j=1}^n j^{-2}+(n^{-1}R_{n\rho})^2\bigg(\frac{1}{n}\sum_{j=1}^N \lambda_i^{\beta-\alpha}\bigg)\bigg)=o(1).$$

Thus, $U_{2n}(P)=o(1)$ in this case.

*Subcase* 7.2. Without loss of generality, we assume $|\sum_{a_{st}\in P_1} a_{st}^{ij}|=\lambda_j+\lambda_i$, $|\sum_{a_{st}\in P_2} a_{st}^{ij}|=\lambda_j-\lambda_i$. Applying Lemma 1(ii) with Lemma 3(i) gives $\prod_{k=3}^4 |\mathrm{cum}(d_{nc}(a_{st}^{ij}):a_{st}\in P_k)|\leq Cn^2(\lambda_i\lambda_j)^{-\alpha}$. Using this, along with Lemma 2, Lemma 3(i) and (20) for a fixed $0<\rho<\pi/2$, we may bound $U_{2n}(P)\leq u_{5n}(\rho)+u_{6n}(\rho)$, where

$$u_{5n}=\frac{C(\rho)}{n^4}\sum_{A^\rho}\lambda_i^{\beta-\alpha}[L_{n0}(\lambda_i+\lambda_j)+R_{n\rho}]\bigg[\frac{n}{j-i}+R_{n\rho}\bigg]$$

and $u_{6n}(\rho)=Cn^{-4}\sum_{A_\rho}(\lambda_i\lambda_j)^{2\beta-\alpha}\lambda_i^{-2\alpha}[\lambda_j^{-1}+L_{n1}(\lambda_i+\lambda_j)][\lambda_j^{-1}+L_{n1}(\lambda_i-\lambda_j)]$. By Lemma 1(ii), we can show $u_{6n}(\rho)\leq Cu_{2n}(\rho)$ while, for large $n$, $u_{5n}(\rho)$ is bounded by

$$\frac{C(\rho)}{n^2}\bigg(\sum_{A^{n/2}}\frac{\lambda_i^{\beta-\alpha}}{n-j-i}\bigg[\frac{1}{j-i}+\frac{R_{n\rho}}{n}\bigg]+\sum_{A_{n/2}}\frac{\lambda_i^{\beta-\alpha}}{i+j}\bigg[\frac{1}{j-i}+\frac{R_{n\rho}}{n}\bigg]$$

$$+R_{n\rho}\sum_{j=1}^N \lambda_j^{\beta-\alpha}\bigg)$$

$$\leq C(\rho)\bigg(n^{-1+\max\{0,\alpha-\beta\}}\sum_{j=1}^n j^{-2}+\frac{R_{n\rho}}{n}\bigg(\frac{1}{n}\sum_{j=1}^N \lambda_i^{\beta-\alpha}\bigg)\bigg)=o(1),$$

using $n-j-i\geq j-i$ above. Hence, $U_{2n}(P)=o(1)$.

*Subcases* 7.3–7.4. For Subcase 7.3, there exists some $k$ such that $|\sum_{a_{st}\in P_k} a_{st}^{ij}|=\lambda_j-\lambda_i$ or $\lambda_j+\lambda_i$. Note that $|\mathrm{cum}(d_{nc}(a_{st}^{ij}):a_{st}\in P_{k_1})|\leq$

$Cn\lambda_i^{-\alpha}$, $|\text{cum}(d_{nc}(a_{st}^{ij}) : a_{st} \in P_{k_2})| \leq Cn\lambda_j^{-\alpha}$ by Lemma 3(i). For Subcase 7.4, the possible partitions $P$ are as follows $|\sum_{a_{st} \in P_k} a_{st}^{ij}| = \lambda_j + \lambda_i$ for each $k$ (or $\lambda_j - \lambda_i$ for each $k$) or there exist $k_1, k_2$ where $|\sum_{a_{st} \in P_{k_1}} a_{st}^{ij}| = \lambda_j + \lambda_i$, $|\sum_{a_{st} \in P_{k_2}} a_{st}^{ij}| = \lambda_j - \lambda_i$. Arguments as in Subcases 7.1–7.2 can show that (34) holds; see Lemma 13 of [32]. $\square$

**Acknowledgments.** The authors wish to thank two referees and an Associate Editor for many constructive comments and suggestions that significantly improved an earlier draft of this paper.

DEPARTMENT OF STATISTICS
IOWA STATE UNIVERSITY
AMES, IOWA 50011
USA
E-MAIL: dnordman@iastate.edu
　　　　snlahiri@iastate.edu